%% file: main.tex
\documentclass[12pt]{article}

\usepackage[english]{babel}
\usepackage[utf8x]{inputenc}
\usepackage{pdfpages}
\usepackage{amsfonts,dsfont}
\usepackage{amsmath}
\usepackage{amssymb}
\usepackage{amsthm} %
\usepackage{bm}
\usepackage{algorithm2e}
\usepackage{todonotes}
\usepackage[all]{xy} 
\usepackage{hyperref}
\usepackage[capitalise]{cleveref}
\usepackage{autonum}

\makeatletter \let\cl@part\relax \makeatother

\makeatletter
\newcommand{\restore@Environment}[1]{%
  \AtBeginDocument{%
    \csletcs{#1*}{#1}%
    \csletcs{end#1*}{end#1}%
  }%
}
\forcsvlist\restore@Environment{alignat,equation,gather,multline,flalign,align}
\makeatother

\usepackage[top=2cm, bottom=2cm, left=2cm, right=2cm]{geometry}
\usepackage{xcolor}
\usepackage{tikz}
\usetikzlibrary{calc}
\usetikzlibrary{intersections}
\tikzset{offset/.style={to path={%
    -- ($(\tikztostart)!#1cm!(\tikztotarget)$)}},
         offset/.default=1}
\tikzset{>=latex}

\usepackage{subcaption}
\usepackage{tikz-3dplot}

\newtheorem{hypo}{Assumption}
\Crefname{hypo}{Assumption}{Assumptions}

\newtheorem{theorem}{Theorem}
\newtheorem{lemma}[theorem]{Lemma}
\newtheorem{remark}[theorem]{Remark}
\newtheorem{cor}[theorem]{Corollary}
\newtheorem{proposition}[theorem]{Proposition}
\newtheorem{definition}[theorem]{Definition}

\theoremstyle{remark}
\newtheorem{example}{Example}

\newcommand{\qedb}{}

\title{Exact quantization of \\
  multistage stochastic linear problems
}
\author{Maël Forcier\footnote{CERMICS, École des Ponts, mael.forcier@enpc.fr}, Stéphane Gaubert\footnote{INRIA, CMAP, École polytechnique, IP Paris, CNRS, stephane.gaubert@inria.fr}, Vincent Leclère\footnote{CERMICS, École des Ponts, vincent.leclere@enpc.fr}}

\input{./commandes}
\def\scalenormequiv{0.6}

\def\scalechcmplx{5}
\def\scaledynprogchcmplx{3}

\def\scalecoupP{2.6}
\def\scalefibers{1.2}
\def\scalefibersfans{0.3}

\begin{document}

\include{main_in}
\bibliographystyle{alpha}
\bibliography{biblio}

\end{document}

%% file: commandes.tex
\newcommand{\ectg}{V} %
\newcommand{\ctg}{\hat{V}} %
\newcommand{\secfan}[1]{\Sigma\operatorname{-fan}(#1)}
\newcommand{\compfansupp}{\cR} %

\newcommand{\nhalfline}[1]{N_{#1}} %

\newcommand{\mynorm}[1]{\|#1\|}

\newcommand{\BB}{{\mathbb B}}

\newcommand{\EE}{{\mathbb E}} %
\newcommand{\II}{\mathbb{I}}

\newcommand{\PP}{{\mathbb P}} %
\newcommand{\QQ}{{\mathbb Q}} %
\newcommand{\RR}{{\mathbb R}} %
\def\SS{\mathbb S}

\newcommand{\ZZ}{{\mathbb Z}} %

\newcommand{\cA}{\mathcal{A}}

\newcommand{\cC}{\mathcal{C}}

\newcommand{\cF}{\mathcal{F}}

\newcommand{\cL}{\mathcal{L}}

\newcommand{\cN}{\mathcal{N}}
\newcommand{\cP}{\mathcal{P}}
\newcommand{\cQ}{\mathcal{Q}}
\newcommand{\cR}{\mathcal{R}}
\newcommand{\cS}{\mathcal{S}}
\newcommand{\cT}{\mathcal{T}}

\newcommand{\cV}{\mathcal{V}}

\newcommand{\cFlow}{\cF_{\mathrm{low}}}

\newcommand{\prj}{\pi}
\newcommand{\prjx}{\pi}

\newcommand{\xyzprjxy}{\pi^{x,y,z}_{x,y}}
\newcommand{\xyzprjx}{\pi^{x,y,z}_{x}}
\newcommand{\xyprjx}{\pi^{x,y}_x}

\newcommand{\xyprjy}{\pi^{x,y}_{y}}
\newcommand{\yzprjy}{\pi^{y,z}_{y}}

\newcommand{\indi}[1]{\mathds{1}_{#1}}
\newcommand{\tropindi}[1]{\II_{#1}}
\newcommand{\centr}[1]{\operatorname{Ctr}\big(#1\big)}
\newcommand{\sphcentr}[1]{\operatorname{SpCtr}\big(#1\big)}
\newcommand{\ang}[1]{\operatorname{Ang}\big(#1\big)}
\newcommand{\Ecind}[1]{\besp{\indi{\va{c} \in #1}\va{c}^\top}}

\newcommand{\chcmplx}[2]{\cC(#1,#2)} %
\newcommand{\chcmplxmax}[2]{\collecmax{\cC}(#1,#2)}

\newcommand{\cansph}[1]{\SS_{#1}} %
\newcommand{\canball}[1]{\BB_{#1}} %

\newcommand{\dimQ}{d}

\newcommand{\expvalf}[2]{\Phi_{#1}(#2)}
\newcommand{\expvalpar}{\theta}
\newcommand{\expval}[1]{\expvalf{#1}{\expvalpar}}
\newcommand{\enclgth}[1]{\langle#1\rangle}
\newcommand{\epiQ}{\epi(Q)}
\newcommand{\afregQ}{\cQ}

\newcommand{\aff}{\operatorname{Aff}}

\newcommand{\card}[1]{\sharp #1}
\newcommand{\collecmax}[1]{#1^{\max}}
\newcommand{\cone}{\operatorname{Cone}}

\newcommand{\dom}{\operatorname{dom}}
\newcommand{\epi}{\operatorname{epi}}
\newcommand{\ie}{\emph{i.e.}}
\newcommand{\im}{\operatorname{Im}}
\renewcommand{\ker}{\operatorname{Ker}}

\newcommand{\rc}{\operatorname{rc}}
\newcommand{\relint}{\operatorname{ri}}
\newcommand{\st}{\text{s.t. }}
\newcommand{\supp}{\operatorname{supp}}

\newcommand{\ray}{\operatorname{Ray}}

\newcommand{\vol}{\operatorname{Vol}}
\newcommand{\vrtx}{\operatorname{Vert}}

\newcommand{\spacevaint}{L^1(\Omega,\cA,\PP;\RR^m)}
\newcommand{\spacevaintt}{L^1(\Omega,\cA,\PP;\RR^{n_t})}
\newcommand{\rctg}{R}
\newcommand{\rafreg}{\cR}

\newcommand{\faceleq}{\vartriangleleft}

\newcommand{\refleq}{\preceq}

\renewcommand{\geq}{\geqslant}
\renewcommand{\leq}{\leqslant}
\renewcommand{\preceq}{\preccurlyeq}
\renewcommand{\succeq}{\succcurlyeq}

{\normalsize}

{\normalsize}

{\normalsize}

\newcommand{\eqfinp}{\; .}

\renewcommand{\tilde}{\widetilde}

\newcommand{\defegal}{:=}
\newcommand{\bp}[1]{\big(#1\big)}                           %
\newcommand{\Bp}[1]{\Big(#1\Big)}                           %
\newcommand{\bc}[1]{\big[#1\big]}                           %
\newcommand{\Bc}[1]{\Big[#1\Big]}                           %
\newcommand{\ba}[1]{\big\{#1\big\}}                         %
\newcommand{\Ba}[1]{\Big\{#1\Big\}}                         %

\newcommand{\horizon}{T}

\newcommand{\argmin}{\mathop{\arg\min}} 
\newcommand{\argmax}{\mathop{\mathrm{argmax}}}

\newcommand{\prbt}{\mathbb{P}}                              %
\newcommand{\mathscr}{\EuScript}

\newcommand{\va}[1]{\bm{#1}}

\newcommand{\bprob}[1]{\prbt\bc{#1}}                  %

\newcommand{\espe}{\mathbb{E}}                              %

\newcommand{\besp}[2][]{\espe_{#1}\bc{#2}}            %
\newcommand{\Besp}[2][]{\espe_{#1}\Bc{#2}}            %

%% file: main_in.tex
\maketitle

\input{sections/abstract}

\section{Introduction}
\label{sec:intro}

\input{sections/intro}
\section{Polyhedral tools}
\label{sec:polyhedral_tools}

\input{sections/polyhedral_tools.tex}
\section{Exact quantization of the $2$-stage problem}
\label{sec:polyhedrality}

\input{sections/preserving_polyhedrality}

\subsection{Alternative approach in terms of regular subdivisions}
\label{sec:basis_decomposition_theorem}

\input{sections/basis_decomposition}

\section{Exact quantization of the multistage problem}
\label{sec:extension}
\input{sections/extensions}

\section{Computing the quantized costs and probabilities}

\label{sec:usual_distr}

\input{sections/formulas}

\subsection{An illustrative example}
\label{sec:analytical_ex}

\input{sections/example}

\section{Complexity}
\label{sec:complexity}

\input{sections/complexity}

%% file: sections/abstract.tex
\begin{abstract}
  We show that the multistage linear problem (MSLP) with an arbitrary
  cost distribution is equivalent to a MSLP on a finite scenario tree.
We establish this exact quantization result by analyzing the polyhedral structure of MSLPs. In particular,  
we show that the expected cost-to-go functions are polyhedral and affine
on the cells of a chamber complex, which is independent
of the cost distribution.
This leads to new complexity results, showing that MSLP
is fixed-parameter tractable. 
\end{abstract}

%% file: sections/intro.tex
Stochastic programming is a powerful modeling paradigm for optimization under uncertainty that has found many applications in energy, logistics or finance (see e.g.\ \cite{wallace2005applications}).
Multistage linear stochastic programs (MLSP) constitute an important class of stochastic programs.
They have been thoroughly studied, see e.g.\ \cite{%
birge2011introduction,prekopa2013stochastic}. %
One reason for this interest is the availability of efficient linear solvers and the use of dedicated algorithms leveraging the special structure of linear stochastic programs (\cite{van1969shaped,birge1985decomposition}).

In this paper,
we show that every MSLP with general cost distribution is equivalent to an MSLP with finite distribution.
This leads to explicit representations of their value functions and to
new complexity results.

\subsection{Multistage stochastic linear programming}

Let $(\Omega,\cA,\PP)$ be a probability space. 
Given a sequence of independent random variables $\va c_t  \in L_1(\Omega,\cA,\PP;\RR^{n_t})$ and 
$\va \xi_t = (\va A_t, \va B_t, \va h_t)$, indexed by $t \in [\horizon]\defegal \{1,\dots,\horizon\}$, we consider the MSLP given by 
\begin{equation}
\label{eq:def_mslp}
\tag{MSLP}
\begin{aligned}
\min_{(\va x_t)_{t\in [\horizon]}}&  \quad  c_1^\top x_1 + \besp{\sum_{t=2}^{\horizon} \va c_t^\top \va x_t  }\\
 \st  
 &\quad A_1 x_1 \leq b_1 \\
 &   \quad \va A_t \va x_{t}+  \va B_t \va x_{t-1} \leq  \va b_t   \quad \text{a.s.} & \forall  t \in \{2,\cdots,\horizon \}\\
 & \quad \va x_t \in L_\infty(\Omega,\cA,\PP;\RR^{n_t}) & \forall  t \in \{2,\cdots,\horizon \}\\
  & \quad \va x_t \preceq \cF_t & \forall  t \in \{2,\cdots,\horizon \} 
\end{aligned}
\end{equation}
where $\va x_1\equiv x_1$, $\va A_1 \equiv A_1$ and $\va b_1 \equiv b_1$ are deterministic and $\cF_t$ is the $\sigma$-algebra generated  by $(\va c_2,\va \xi_2,\cdots, \va c_{t}, \va \xi_{t})$.
The last constraint, known as non-anticipativity, means that $\va x_t$ is measurable with respect to $\cF_t$.

Most results for MSLP with continuous distributions rely on discretizing the distributions. The Sample Average Approximation (SAA) method (see e.g.\ \cite[Chap. 5]{shapiro2014lectures}) samples the costs and constraints.
It relies on probabilistic results based on a uniform law of large number to give statistical guarantees. 
Obtaining a good approximation requires a large number of scenarios.
In order to alleviate the computations, we can use
 scenario reduction techniques (see \cite{dupacova2003scenario,heitsch2003scenario}).
Latin Hypercube Sampling and variance reduction methods are also used to produce scenarios. 
Finally, one generates heuristically ``good'' scenarios,
representing the underlying distribution (see \cite{KautWall}).
Alternatively, we can leverage the structure of the problem to produce finite scenario trees (see \cite{kuhn2006generalized,maggioni2014bounds,Maggioni2018}) that yields bounds for the value of the true optimization problem.
In each of these approaches, one solves 
an approximate version of the stochastic program, with or without statistical guarantee.

\subsection{The exact quantization problem}
Here, we aim at solving exactly the original problem,
by finding an equivalent formulation with discrete distributions. This notion
of equivalent formulation is best understood through the dynamic programming approach
of MSLP. We define the \emph{cost-to-go} function $\ectg_t$ inductively as follows.
We set %
$\ectg_{\horizon+1} \equiv 0$ and for all $t\in \{2,\cdots,\horizon \}$:
\begin{equation}
\label{eq:def_ectg_multi}
\begin{aligned}
\ectg_t(x_{t-1}) &\defegal \besp{\ctg_t(x_{t-1},\va c_t,\va \xi_t)}\\
    \ctg_t(x_{t-1},c_t,\xi_t)&\defegal \min_{x_t \in \RR^{n_t}} \quad  c_t^\top x_t + \ectg_{t+1}(x_t)\\
    & \quad \quad \st \quad A_t x_t +  B_tx_{t-1} \leq  h_t 
\end{aligned}
\end{equation}
where $x_{t-1} \in \RR^{n_{t-1}}$, $c_t \in \RR^{n_t}$ and $\xi_t = (A_t,B_t,b_t) \in \RR^{q_t \times n_{t-1}} \times \RR^{q_t \times n_t} \times \RR^{q_t}$.

We choose to distinguish the random cost $\va c_t$ from the noise $\va \xi_t$ affecting the constraints. 
Indeed our results require $\va \xi_t$ to be finitely supported (see \cref{ex:stoch_h,ex:stoch_T}) while $\va c_t$ can have a continuous distribution.
This separation does not preclude correlation between $\va c_t$ and $\va \xi_t$.
However, we require $\{(\va c_t, \va \xi_t)\}_{t\in [\horizon]}$ to be a sequence of independent random variables to leverage Dynamic Programming, even though some results can be extended to dependent $(\va \xi_t)_{t\in[\horizon]}$.

We say that an MSLP admits an \emph{exact quantization} if there exists a finitely supported $(\check{\va c}_t,\check{\va\xi}_t)_{t\in [\horizon]}$ that yields the same
expected cost-to-go functions, $(\ectg_t)_{t\in [\horizon]}$.
In particular the MSLP is equivalent to a problem on a finite scenario tree.

An obvious necessary condition for exact quantization
is that the value function $\ectg_t$ be a polyhedral function,
meaning that it takes value in $\RR \cup \{+\infty\}$ and its epigraph is a (possibly empty) polyhedron.
Indeed, for each $(c,\xi) \in \supp (\va c,\va \xi)$, 
$Q^{c,\xi} : (x,y)\to  \; c^\top y +\tropindi{A x+ B y\leq h}$ is polyhedral.
Thus, $\ctg(\cdot,c,\xi):=\min_{y\in \RR^m} Q^{c,\xi}(\cdot,y)$ is polyhedral as $\epi \ctg(\cdot,c,\xi)$ is a projection of $\epi Q^{c,\xi}$ (see \cite{jones2008polyhedral,shapiro2014lectures}).
Hence, the following examples show that if the constraints have non-discrete distributions, there is no hope to have an exact quantization theorem.
We shall see, however, that this is the case
without restrictions on the cost distribution.

\begin{example}[Stochastic $\va B$]
\label{ex:stoch_T}
Here, and in the next example, $\va u$ denotes a uniform random variable on $[0,1]$.
\begin{equation}
\ectg(x)
    =  \EE\left[
    \begin{aligned}\min_{y \in \RR^m} \quad & y  \\
    \st \quad & \va{u}x\leq y\\
    & 1\leq y
    \end{aligned}\right]
    = \besp{\max(\va{u}x,1)}
    =  \begin{cases}
    1 & \text{if } x\leq 1\\
    \frac{x}{2}+\frac{1}{2x} & \text{if } x\geq 1
    \end{cases}
\end{equation}
\end{example}

\begin{example}[Stochastic $\va b$]
\label{ex:stoch_h}
\begin{equation}
\ectg(x)
    =  \EE\left[
    \begin{aligned}\min_{y \in \RR^m} \quad & y  \\
    \st \quad& \va{u}\leq y\\
    & x\leq y
    \end{aligned}\right]
     = \besp{\max(x,\va{u})}
    =  \begin{cases}
    \frac{1}{2} & \text{if } x\leq 0\\
    \frac{x^2+1}{2} & \text{if } x\in [0,1]\\
    x & \text{if } x\geq 1 \enspace . 
    \end{cases}
\end{equation}

\end{example}

\subsection{Contribution}

We rely on a geometric approach, which
enlightens the polyhedral structure of MSLP.
We first establish exact quantization results in 
the two-stage case showing that there exists an optimal recourse affine on each cell of a polyhedral complex
which is precisely the chamber complex \cite{billera1992fiber,rambau1996projections}, a fundamental object in combinatorial geometry.
A chamber complex is 
defined as the common refinement of the projections
of faces of a polyhedron.
In particular, \Cref{thm:quantization_2stage}
provides an explicit exact quantization, in which the quantized
probabilities and costs are attached to the cones
of a polyhedral fan $\cN$
(we refer the reader to \cite{de2010triangulations,ziegler2012lectures,grunbaum2013convex,fukuda2016lecture}
for background on polyhedral complexes and fans).
On each cone $N \in \cN$, we replace the distribution of $\va c \indi{ \relint{N}}$ by a Dirac distribution concentrated on
the expected value $\check{c}_N=\besp{\va c | \va c \in \relint{N}}$, and an associated weight $\check{p}_N=\bprob{\va c \in \relint{N}}$. 
Further, $\cN$ is \emph{universal} in the sense that it does not depends on the distribution of $\va c$.  

In order to extend this result to the multistage case
we establish in \Cref{lem:DP_chamber_complex} a Dynamic Programming type equation in the space of polyhedral complexes. 
Then we show an exact quantization result in \Cref{thm:multistage_quantization}.
Again, this quantization is \emph{universal} in the cost distribution. 

We apply this polyhedral approach to obtain fixed 
parameters polynomial time complexity results considering
both the exact computation problem and the approximation problem. 
For distributions that are uniform on polytopes or exponential,
we show the MLSP can be solved in a time that is polynomial
provided that the horizon $\horizon$ and the dimensions $n_2,\dots,n_\horizon$ of the successive recourses are fixed.
The proof relies on the theory of linear programming with oracles~\cite{grotschel2012geometric} as well as on upper 
bound theorems of McMullen~\cite{mcmullen1970maximum} and Stanley~\cite{stanley1975upper} concerning the number of vertices and the size of a triangulation of a polyhedron.
We obtain similar results for the approximation problem. 
Then the distribution cost can be essentially arbitrary: we only assume that it is given implicitely through an appropriate oracle.
This applies in particular to distributions with a smooth density with respect to Lebesgue measure.

In summary, our main contributions are the following:
\begin{enumerate}
    \item MSLP with arbitrary cost distribution and finitely supported constraints admit an exact quantization result, \ie\ are equivalent to MSLP with discrete cost distribution;
    \item the cost-to-go functions of such MSLP are polyhedral and
        affine on the cells of a universal polyhedral complex (\ie\  independent
        of the cost distribution);
    \item exact formulas for quantized costs and probabilities in the case of exponentially or uniformly distributed costs on a polytope; 
    \item fixed-parameter polynomial time tractability results
        for 2SLP and MSLP.
\end{enumerate}

\subsection{Comparison with related work}

A combinatorial approach of deterministic parametric linear programming  was developped by Walkup and Wets \cite{walkup1969lifting} see also \cite{sturmfels1997variation} for a more recent discussion.
Their basis decomposition theorem describes how the value 
of a linear program in standard form varies with respect to 
the cost and to the right-hand side of the constraints. 
In the two stage case, we can see the collection of rows of $A$ as a vector configuration, 
and the right-hand side of the recourse problem $b-Bx$ as a height function which determines a regular subdivision of this configuration.
The space of regular subdivision is represented by the so called secondary fan \cite{de2010triangulations}.
We may apply this theorem to the dual problem of the recourse problem to deduce that the expected cost-to-go function is affine on each cell of an affine section of the secondary fan.
This affine section can be shown to coincide with the chamber complex used here. 
However, the basis decomposition theorem cannot be applied to the extensive form of a \emph{multistage} problem. 
In particular nonanticipativity constraints cannot be tackled in this way.
Thus, we choose to develop an approach through chamber complexes as it is more direct, 
allowing us to obtain also a result in the multistage case. 
The comparison with the approach of Walkup and Wets is further discussed in \Cref{sec:basis_decomposition_theorem}.

The complexity of stochastic programming has been extentively studied. 
Dyer and Stougie~\cite{dyer2006computational} proved that 2 stage stochastic programming is $\sharp P$-hard in the discrete case,
 by reducing the problem of graph reliability to the discrete distribution case.
 They stated that the computation of the volume of a polytope can be reduced to the continuous distribution case,
a result which was subsequently proved in~\cite{hanasusanto2016comment}.
Computing the volume of a polytope, as well as graph reliability,
is $\sharp P$-complete.
Hanasusanto, Kuhn and Wiesemann (\emph{ibid}) showed that computing an approximate solution to the 
2-stage linear programming (2SLP) with continuous distribution with a sufficiently high accuracy is also $\sharp P$-hard.
Other papers \cite{shapiro2005complexity} studied the complexity of 2-stage linear programming 2SLP and MSLP. 
Most complexity results there are hardness results.
In contrast, we prove that 2SLP and MSLP are fixed parameter tractable.

Finally,
    Lan \cite{lan2020complexity} and Zhang and Sun \cite{zhang2019stochastic} independently analysed the complexity of Stochastic Dual Dynamic Programming (SDDP). 
    It follows from their results that
    finitely supported MSLP can be solved approximately in pseudo-polynomial time in the error approximation $\varepsilon$ when all the dimensions and the horizon are fixed.
    In other words the complexity of these SDDP methods is polynomially bounded in $1/\varepsilon$.
In contrast, our approach shows that MSLP can be solved approximately in polynomial time in $\log(1/\varepsilon)$, when $\horizon$, $n_2, \cdots, n_{\horizon}$ are fixed.
In particular, the first state dimension is not fixed.
Moreover, we obtain polynomial complexity bounds in the exact (Turing) model of computation for appropriate classes of distributions.
Note that in the approach presented here, contrary to SDDP like methods, we do not rely on statistical sampling and the value functions are computed exactly in one pass only. 
However, the objective of SDDP is to obtain quickly an approximate solution whereas our approach computes exactly all the supporting hyperplanes.

\subsection{Structure of the paper}
We recall, in~\Cref{sec:polyhedral_tools}, notions from the theory of polyhedra:
\emph{polyhedral complexes}, \emph{normal fans} and \emph{chamber complexes}.
    In \Cref{sec:polyhedrality} we establish the exact quantization result for 2-stage stochastic linear programming.
    In \Cref{sec:extension}, we show that  chamber complexes can be propagated through dynamic programming,
    leading to the exact quantization result for the MSLP.
We show in~\Cref{sec:usual_distr} how the quantized probabilities and cost can be computed for appropriate distributions.
Finally, in~\Cref{sec:complexity}, we draw the consequences of our results in terms of computational complexity.

\subsection{Notation}
\label{sec:notation}

As a general guideline $\va{bold}$ letters denote random variables, normal scripts their realisation.
Capital letters denote matrices or sets, calligraphic (e.g.\ $\cN$) denote
collections of sets.
The indicator function $\tropindi{\mathtt{P}}$ (resp.\ $\indi{\mathtt{P}}$) takes value $0$ (resp.\ $1$) if $P$ is true and $+\infty$ (resp.\ $0$) otherwise.
We set $[k]:=\{1,\dots,k\}$, and we denote by $\card E$ the cardinal of a set $E$.
We denote by $\cone(A):=A\RR_+^n$ the cone hull of the columns of $A$.
$x \leq y$ is the standard partial order, given by $\forall i, x_i\leq y_i$. 
$F \faceleq G$ if $F$ is a subface of $G$.
    $\cP \refleq \cQ$ if $\cP$ is a refinement of the polyhedral complex $\cQ$.
$\supp \cC := \bigcup_{C \in \cC} E$ is the support of a collection of sets $\cC$, $\collecmax{\cC}$ : the sets of maximal elements of a collection of sets $\cC$.
$\rc(P)$ is the recession cone of a polyhedron $P$.
For a polyhedron $P$, we denote $\cF(P)$ its faces, $\vrtx(P)$ its vertices and $\ray(P)$ a set with vectors each representing one extreme rays (for example the normalized extreme rays).
$P^\psi$ is the face of $P$ given by $\argmin_{x \in P} \psi^\top x$.
$N_P(x)$ is the normal cone of $P$ at $x$, and $\cN(P)$ the normal fan of $P$.

%% file: sections/polyhedral_tools.tex
Our proofs rely on the notions of normal fan and chamber complex of a polyhedron recalled here.
These polyhedral objects reveal the geometrical structure of MSLP.
Both the normal fan and the chamber complex are special polyhedral complexes.

\subsection{Polyhedral complexes}

\emph{Polyhedral complexes} are collections of polyhedra 
satisfying some combinatorial and geometrical properties.
In particular the relative interiors of the elements of a polyhedral complex (without the empty set) form a partition of their union.
We refer to \cite{de2010triangulations} for a complete introduction to polyhedral complexes and triangulations.

\begin{definition}[Polyhedral complex]
A finite collection of polyhedra $\cC$ is a \emph{polyhedral complex} 
if it satisfies 
i) if $P \in \cC$ and $F$ is a non-empty\footnote{For some authors, a polyhedral complex must contain the empty set. We do not make this requirement.} face of $P$ then $F \in \cC$ and 
ii) if $P$ and $Q$ are in $\cC$, then $P\cap Q$ is a (possibly empty) face of $P$.

We denote by $\supp \cC:=\bigcup_{P \in \cC} P$ the support of a polyhedral complex.
Further, if all the elements of $\cC$ are polytopes (resp. cones, simplices, simplicial cones), 
we say that $\cC$ is a \emph{polytopal complex} (resp. a \emph{fan}, a \emph{simplicial complex}, a \emph{simplicial fan}).
\end{definition}

We recall that a \emph{simplex} of dimension $d$ is the convex hull of $d+1$ affinely independent point and that a \emph{simplicial cone} of dimension $d$ is the conical hull of $d$ linearly independent vectors. 
\begin{proposition}
\label{prop:relint_partition}
For any polyhedral complex $\cC$, the relative interiors of its elements (without the empty set) form a partition of its support: $\supp(\cC) = \bigsqcup_{P \in \cC} \relint(P)$.
\end{proposition}

For example, the set of faces $\cF(P)$ of a polyhedron $P$ is a polyhedral complex.

\begin{definition}[Refinements and triangulation]
\label{def:refinement}
Let $\cC$ and $\cR$ be two polyhedral complexes, we say that $\cR$ is a \emph{refinement} of $\cC$, 
denoted $\cR \refleq \cC$, if $\supp \cR=\supp \cC$ and for every cell $R \in \cR$ there exists a cell $C \in \cC$ containing $R$: $R \subset C$.

Note that $\refleq$ defines a partial order on the space of polyhedral complexes, and the \emph{meet}
associated with this order is given by
the \emph{common refinement} of two polyhedral complexes $\cC$ and $\cC^\prime$ 
defined as the polyhedral complex of the intersections of cells of $\cC$ and $\cC^\prime$:
\begin{equation}
\cC \wedge \cC^\prime := \{ R\cap R^\prime \, | \, R \in \cC, R^\prime \in \cC^\prime\}
\end{equation}

	A \emph{triangulation $\cT$ of a polytope} $Q$ is a refinement of $\cF(Q)$ such that the cells of dimension $0$ of $\cT$ are the vertices of $Q$ and $\cT$ is a simplicial complex.
A \emph{triangulation $\cT$ of a cone} $K$ is a refinement of $\cF(K)$ such that the cells of dimension $1$ of $\cT$ are the rays of $K$ and $\cT$ is a simplicial fan.

\end{definition}

\subsection{Normal fan}

The normal fan is the collection of the normal cones of all faces of a polyhedron.
See \cite{lu2008normal} for a review of normal fan properties. 

Recall that the \emph{normal cone} of a convex set $C \subset \RR^m$ at the point $x$ is the set
$N_C(x):=\{\alpha \in \RR^m\; |\; \forall y \in C, \; \alpha^\top (y-x)\leq 0\}$.
More generally, for a set $E \subset C$, 
$N_C(E):=\bigcap_{x \in E}N_C(x)$.%

\begin{figure}[ht]
\input{tikz/ex_normally_equivalent}
\end{figure}

\begin{definition}[Normal Fan]
The \emph{normal fan}\footnote{Sometimes called \emph{outer} normal cones and fan, as opposed to \emph{inner} cones obtained either by inverting the inequality in the definition of the normal cone or by taking the opposite cones respect to the origin.
} 
of a convex set $C$ is the collection of polyhedral cones
\begin{equation}
\cN(C):=\{ N_C(x) \; |\; x \in C\}
\end{equation}
We say that two convex sets $C$ and $C^\prime$ are \emph{normally equivalent} if they have the same normal fan : $\cN(C)=\cN(C^\prime)$.
\end{definition}

Recall that the \emph{polar} of a convex set $C$ is the set
$
C^\circ:=\{\alpha\; |\; \forall x \in C, \;\alpha^\top x\leq 0 \} = N_C(0)
$
and the {\em recession cone} of a convex set $C$ is given by 
$\rc(C):=\{ r \in C \;|\; \forall \mu \in \RR_+, \; \forall x \in c, \; x + \mu r \in C \}$.
In particular, for a polyhedron, the recession cone 
and its polar are given by 
\begin{equation}
\rc\bp{\{x\;|\; Ax\leq b\}}=\{x\;|\; Ax\leq 0\} 
\qquad 
\rc\bp{\{x\;|\; Ax\leq b\}}^\circ = \cone(A^\top)
\eqfinp
\label{eq:rc_polyhedron}
\end{equation}

\begin{proposition}[Basic properties of normal fans (see e.g.\ \cite{lu2008normal})]
\label{prop:fan_basic}

If $P$ is a polyhedron, the normal fan $\cN(P)$ is a
\emph{finite} collection of polyhedral cones (and in particular
a \emph{polyhedral complex}). 
Further, the support of $\cN(P)$ can be expressed geometrically as the polar of the recession cone of $P$, \ie
\begin{equation}
\supp \cN(P) = \big(\rc(P)\big)^\circ
\label{eq:supp_fan_rc}
\end{equation}

\end{proposition}

\subsection{Chamber complex}

The affine regions of the cost-to-go function will correspond to cells of a chamber complex.
Projections of polyhedra, fibers and chambers complexes are studied in \cite{billera1992fiber,rambau1996projections,rambau1996polyhedral}.

\begin{definition}[Chamber complex]
\label{def:chmbrcplx}
Let $P \subset \RR^n$ be a polyhedron and $\prj$ a linear projection defined on $\RR^n$.
For $x \in \prj(P)$ we define the \emph{chamber of $x$ for $P$ along $\prj$} as
\begin{equation}
\sigma_{P,\prj}(x):= \bigcap_{F \in \cF(P)\, \st\, x \in \prj(F)} \prj(F) .
\end{equation}
The \emph{chamber complex} $\chcmplx{P}{\prj}$ of $P$ along $\prj$ is  defined as the (finite) collection of chambers, \ie
\begin{equation}
\chcmplx{P}{\prj}:=\{ \sigma_{P,\prj}(x) \; |\; x \in \prj(P)\} \, .
\end{equation}
Further $\chcmplx{P}{\prj}$ is a polyhedral complex such that $ \supp \chcmplx{P}{\prj} = \prj(P)$.
In particular, $\ba{\relint(\sigma) \, | \, \sigma \in \chcmplx{P}{\prj}  }$ is a partition of $\prj(P)$.

More generally, the chamber complex of a polyhedral complex $\cP$ is 
\begin{equation}
\label{eq:chcmplx_cplx_ploy}
\chcmplx{\cP}{\prj}:=\{ \sigma_{\cP,\prj}(x) \; |\; x \in \prj\bp{\supp(\cP)}\} \, .
\end{equation}
with $\sigma_{\cP,\prj}(x):= \bigcap\limits_{F \in \cP\, \st\, x \in \prj(F)} \prj(F)$.

\end{definition}

\begin{figure}[ht]
\input{tikz/ex_chamber_complex}
\label{fig:chamber_complex_ex}
\end{figure}

    \begin{lemma}[Chamber complex monotonicity with respect to refinement order]
    \label{lem:chmbrcplx_refinement}
        Consider two polyhedral complexes of $\RR^d$ and a projection $\prj$.
 If $\cR \refleq \cS$ then $\chcmplx{\cR}{\prj} \refleq \chcmplx{\cS}{\prj}$.
    \end{lemma}
    \begin{proof}

For any $R \in \cR$, there exist $S_R \in \cS$ such that $R\subset S_R$.
Let 
$x\in \supp\chcmplx{\cR}{\prj}=\prj (\supp \cR)=\prj (\supp \cS)=\supp\chcmplx{\cS}{\prj}$
\begin{align*}
\sigma_{\cR,\prj}(x):=&\bigcap_{R \in \cR \, \st x\in \prj(R)}\prj(R)
\subset \bigcap_{R \in \cR \, \st x\in \prj(R)}\prj(S_R) \\
&\subset \bigcap_{S \in \cS \st x\in \prj(S)}\prj(S)
=:\sigma_{\cS,\prj}(x) \in \chcmplx{\cS}{\prj} 
\end{align*}
\qedb
    \end{proof}

Recall that the \emph{fiber} $P_x$ of $P$ along $\prj$ at $x$ is the projection of  $P\cap \prj^{-1}(\{x\})$  on the space $\ker(\prj)$ (see figure~\ref{fig:chamber_complex_ex}). 
An important property of a chamber complex is that all fibers are normally equivalent in each relative interior of cells of the chamber complex. 
More precisely, let $\sigma \in \chcmplx{P}{\prj}$ be a chamber, and $x$ and $x'$ two points in its relative interior, 
then,
$P_x$ and $P_{x^\prime}$ are \emph{normally equivalent}, i.e.\ they have the same normal fan 
$\cN(P_x)=\cN(P_{x^\prime})$,
see~\cite{billera1992fiber}.
Thus we define the normal fan $\cN_\sigma$ above\footnote{The normal fan $\cN_\sigma \subset 2^{\RR^m}$ above $\sigma$ should not be confused with $\cN(\sigma)\subset 2^{\RR^n}$ the normal fan of $\sigma$ which will never appear in this paper.} $\sigma\in \chcmplx{P}{\prj}$ by :
\begin{equation}
\cN_\sigma := \cN(P_x) \quad \text{for an arbitrary } x \in \relint(\sigma)
\end{equation}
The terms {\em parametrized polyhedron}, instead of fibers, and {\em validity domains}, instead of chambers, are also used in the literature \cite{clauss1998parametric,loechner1997parameterized}.

%% file: tikz/ex_normally_equivalent.tex
\def\raycone{1}
\def\lengthcone{1.1}
\def\colorcone{green!40!white}
\def\colorboundcone{green!70!blue}
\def\colorreccone{red}
\def\colorpoly{blue!40!white}
\def\colorbigpoly{blue}
\def\colorcuttingplane{yellow!60!white}
\def\origin{0,0,0}

\def\ax{0}      \def\ay{4}
\def\bx{2}      \def\by{5}
\def\cx{5}      \def\cy{4}
\def\dx{4}      \def\dy{2}
\def\ex{1}      \def\ey{0}

\def\abx{\ax-\bx}      \def\aby{\ay-\by}
\def\bcx{\bx-\cx}      \def\bcy{\by-\cy}
\def\cdx{\cx-\dx}      \def\cdy{\cy-\dy}
\def\dex{\dx-\ex}      \def\dey{\dy-\ey}
\def\eax{\ex-\ax}      \def\eay{\ey-\ay}

\begin{subfigure}{0.45\textwidth}
  \begin{tikzpicture}[scale=\scalenormequiv]
\input{tikz/polytope_fan_macro}
\end{tikzpicture}
\centering
\end{subfigure}
\begin{subfigure}{0.45\textwidth}

\begin{tikzpicture}[scale=\scalenormequiv]
\pgfmathsetmacro\aby{\ay-\by};
\pgfmathsetmacro\bcy{\by-\cy};
\pgfmathsetmacro\cdy{\cy-\dy};
\pgfmathsetmacro\dey{\dy-\ey};
\pgfmathsetmacro\eay{\ey-\ay};

\pgfmathsetmacro\abx{\ax-\bx};
\pgfmathsetmacro\bcx{\bx-\cx};
\pgfmathsetmacro\cdx{\cx-\dx};
\pgfmathsetmacro\dex{\dx-\ex};
\pgfmathsetmacro\eax{\ex-\ax};

\pgfmathsetmacro\ax{\ax+0.9*\abx};
\pgfmathsetmacro\ay{\ay+0.9*\aby};
\pgfmathsetmacro\ex{\ex+0.9*\abx};
\pgfmathsetmacro\ey{\ey+0.9*\aby};

\input{tikz/polytope_fan_macro}
\end{tikzpicture}
\centering
\end{subfigure}%

\caption{Two normally equivalent polytopes $P$ and $P^\prime$ and their normal fan $\cN(P)=\cN(P^\prime)$.}

%% file: tikz/ex_chamber_complex.tex
\def\raycone{1}
\def\lengthcone{1.1}
\def\colorcone{green!40!white}
\def\colorboundcone{green!70!blue}
\def\colorreccone{red}
\def\colorchamber{red}
\def\colorfiber{blue}
\def\colorpoly{green!10!white}
\def\colorbigpoly{blue}
\def\colorcuttingplane{yellow!60!white}
\def\origin{0,0,0}

\def\ax{0.3}    \def\ay{0.6}
\def\bx{0.4}    \def\by{1}
\def\cx{0.8}    \def\cy{1.1}
\def\dx{1.35}   \def\dy{0.8}
\def\ex{1.4}    \def\ey{0.4}
\def\fx{1.2}    \def\fy{0.1}
\def\gx{0.5}    \def\gy{0.3}

\def\xfiber{0.7}

\begin{tikzpicture}[scale=\scalechcmplx]
\coordinate (a) at (\ax,\ay);
\coordinate (b) at (\bx,\by);
\coordinate (c) at (\cx,\cy);
\coordinate (d) at (\dx,\dy);
\coordinate (e) at (\ex,\ey);
\coordinate (f) at (\fx,\fy);
\coordinate (g) at (\gx,\gy);

\filldraw[name path=polyhedron,fill=\colorpoly] (a) -- (b) -- (c) --(d) -- (e) -- (f) -- (g) -- (a);

\draw[thin,dashed] (\ax,0) -- (\ax,\ay) ;
\draw[thin,dashed] (\bx,0) -- (\bx,\by) ;
\draw[thin,dashed] (\cx,0) -- (\cx,\cy) ;
\draw[thin,dashed] (\dx,0) -- (\dx,\dy) ;
\draw[thin,dashed] (\ex,0) -- (\ex,\ey) ;
\draw[thin,dashed] (\fx,0) -- (\fx,\fy) ;
\draw[thin,dashed] (\gx,0) -- (\gx,\gy) ;

\draw[thin,->] (0,0) -- (1.6,0) node[anchor=west] {$\im \prj$};
\draw[thin,->] (0,0) -- (0,1.25) node[anchor=south] {$\ker \prj$};\draw[\colorchamber] (\ax,0)node {$\bullet$} -- (\bx,0) node {$\bullet$} -- (\cx,0) node {$\bullet$} -- (\dx,0) node {$\bullet$} -- (\ex,0) node {$\bullet$}-- (\fx,0) node {$\bullet$}-- (\gx,0) node {$\bullet$};
\draw[\colorchamber] (1,0) node[anchor=south] {$\chcmplx{P}{\prj}$};

\fill[\colorfiber] (\xfiber,0) node[anchor=north]{$x$} circle (0.3pt);

;
\draw[\colorfiber,dashed,name path=verticalx] (\xfiber,0) -- (\xfiber,1.2);
\filldraw[fill=\colorfiber,name intersections={of=polyhedron and verticalx,total=\t}]
    \foreach \s in {1,...,\t}{(intersection-\s) circle (0.3pt) };

\draw[\colorfiber,name intersections={of=polyhedron and verticalx}]
    (intersection-1) -- ($(intersection-1)!0.5!(intersection-2)$) node[anchor=south,rotate=90] {\small$\prj^{-1}(x)\cap P$} --(intersection-2);

\draw[\colorfiber,name intersections={of=polyhedron and verticalx}]
	let \p1=(intersection-1),\p2=(intersection-2) in (0,\y1) node {$\bullet$} -- ($(0,\y1)!0.5!(0,\y2)$) node[anchor=east] {$P_x$}--(0,\y2) node {$\bullet$};

\draw[\colorfiber,dashed,name intersections={of=polyhedron and verticalx}]
	let \p1=(intersection-1) in (\p1) -- (0,\y1) ;

\draw[\colorfiber,dashed,name intersections={of=polyhedron and verticalx}]
	let \p2=(intersection-2) in (0,\y2) -- (\p2);

\draw[violet] (c) node {$\bullet$} -- ($(c)!0.5!(d)$) node[anchor=south west] {$F$}-- (d) node {$\bullet$};

\draw[violet] (\cx,-0.05) node {$\bullet$} -- ($(\cx,-0.05)!0.5!(\dx,-0.05)$) node[anchor=north] {$\prj(F)$}-- (\dx,-0.05) node {$\bullet$};

\draw[green!70!black] (\ax,-0.17) node {$\bullet$} -- ($(\ax,-0.17)!0.5!(\dx,-0.17)$) node[anchor=north] {$\prj(P)$}-- (\ex,-0.17) node {$\bullet$};

\end{tikzpicture}
\centering
\caption{A polytope $P$ in light green, its chamber complex in red on the $x$-axis and a fiber $P_x$ in blue on the $y$-axis, for the orthogonal projection $\prj$ on the horizontal axis.}

%% file: sections/preserving_polyhedrality.tex
Let $(\Omega,\cA,\PP)$ be a probability space, $\va c \in L_1(\Omega,\cA,\PP;\RR^m)$ be an integrable random vector, and suppose $\xi = (T,W,h)$ is deterministic. 
We study the cost-to-go function of the 2-stage stochastic problem, written as
	\begin{equation}
	\label{pb:2nd-stage}
      \begin{aligned}
      \ectg(x) := \EE\left[ \ctg(x,\va c)\right]
     \quad \text{with} \quad \ctg(x,c) := \min_{y \in \RR^m} \quad & c^\top y \\
    \st \quad & Ay+  Bx\leq b 
    \end{aligned}
\end{equation}

	The dual of the latter problem, for given $x$ and $c$, is
	 \begin{align}
			 \max_{\mu \in \RR^q} \label{pb:dual_2nd-stage}
				&\quad (Bx-b)^\top\mu \\
				\st &\quad  A^\top \mu =-c \\
				&\quad \mu\geq 0
	\end{align}

We denote the \emph{coupling constraint polyhedron} of Problem~\eqref{pb:2nd-stage} by
\begin{equation}
P:=\{(x,y) \in \RR^{n+m}\; |\; Ay+Bx\leq b \}
\label{e:def:Twh}
\end{equation}
and $\prjx$ the projection of $\RR^{n} \times \RR^{m}$ onto $\RR^n$ such that $ \prjx(x,y)=x$.

The projection of $P$ is the following polyhedron :
\begin{equation}
\prjx(P)=\{x \in \RR^{n}\; |\; \exists y \in \RR^m, \; Ay+Bx\leq b \}
\label{eq:prj_P}
\end{equation}
and for any $x \in \RR^n$, the \emph{fiber} of $P$ along $\prjx$ is
\begin{equation}
\label{eq:def_Px}
P_x:=\{y \in \RR^{m}\; |\; Ay+Bx\leq b \}
\end{equation}

\subsection{Chamber complexes arising from 2-stage problems}

\label{ssec:finite_scenario}

The following lemma provides an explicit formula for the cost-to-go
function. It shows that an optimal recourse can be chosen
as a function of $c$ that is piecewise constant on
the normal fan of $P_x$.

\begin{lemma}
	Let $x \in \RR^n$ and $c\in \RR^m$,
	\begin{enumerate}
		\item  If $x \notin\prj(P)$, then $\ctg(x,c)=+\infty$;
		\item  If $x \in\prj(P)$ and $-c \notin \cone(A^\top)$, then $\ctg(x,c)=-\infty$;
		\item Suppose now that $x \in\prj(P)$ and $-c \in \cone(A^\top)$. For each

                cone $N\in \cN(P_x)$, let us select in an arbitrary manner a vector $c_N$ in $\relint(-N)$. 
                Then, there exists a vector $y_N(x)$ which achieves the minimum in the expression of $\ctg(x,c_N)$ in~\eqref{pb:2nd-stage}. 
                Further, for any selection of such a $y_N(x)$, we have
                \begin{align}
                \ctg(x,c)= \sum_{N \in \cN(P_x)}  \indi{c \in-\relint N}\; c^\top y_N(x) \label{e-sumformula}
                \enspace .\end{align}

	\end{enumerate}
	\label{lem:value_ctg}
\end{lemma}

\begin{proof}
The first point comes from the definitions of $\prj(P)$ in 
\eqref{eq:prj_P}
and $\ctg(x,c)$ in \eqref{pb:2nd-stage}.
 If $x \in\prj(P)$ and $-c \notin \cone(A^\top)$, then the primal problem \cref{pb:2nd-stage} is feasible and the dual problem 
is \cref{pb:dual_2nd-stage} unfeasible.
Thus, by strong duality, $\ctg(x,c)=-\infty$.

	 By \cref{eq:supp_fan_rc}, we have that $\rc(P_x)^\circ=\supp(\cN(P_x))$.
Further, by \cref{eq:rc_polyhedron} all non empty fibers $P_x$ have the same
recession cone $\{y\;|\; Ay\leq 0\}$ whose polar is $\cone(A^\top)$.

Assume now that $x \in \prj(P)$ and $-c \in \cone(A^\top)=\supp(\cN(P_x))$.
Then, there exists $N \in \cN(P_x)$ such that $-c \in \relint(N)$.
Moreover, for every choice of $c_N \in -\relint(N)$, we have $\argmin_{y \in P_x} c^\top y =\argmin_{y \in P_x} c_N^\top y$, see e.g.\ \cite[Cor.~1(c)]{lu2008normal}.
Moreover, there exists $y_N(x)$ such that $N=N_{P_x}\bp{y_N(x)}$ by definition of a normal cone, thus $y_N(x) \in \argmin_{y \in P_x} c_N^\top y$; in particular, the latter argmin is non empty.
Thus, when $-c\in \relint(N)$, $\ctg(x,c)=c^\top y_N(x)$.

Thanks to the partition property of \cref{prop:relint_partition}, we know that $c$ belongs to the relative interior of precisely one cone in the normal fan of $P_x$,
leading to~\eqref{e-sumformula}.
\end{proof}

Having this property in mind, we make the following assumption:
\begin{hypo}
	The cost $\va c \in \spacevaint$ is integrable with $\va c \in -\cone(A^\top)$ almost surely.
	\label{as:supp_vac}
\end{hypo}

	\begin{theorem}[Quantization of the cost distribution]
\label{thm:quantization_2stage}
Recall that $\chcmplx{P}{\prjx}$ is the chamber complex of the coupling constraint polyhedron $P$ along the projection $\prjx$ on the $x$-space.
Let $x \in \prj(P)$, and $\sigma$ be a
cell of $\chcmplx{P}{\prjx}$ such that $x \in \relint(\sigma)$.

Under \cref{as:supp_vac}, for every refinement $\cR$ of $-\cN_\sigma$, we have:
\begin{equation}
\label{eq:red_scen_det}
\begin{aligned}
      \ectg(x)=  \sum_{R \in \cR} \check p_R \ctg(x,\check c_R)
     \quad \text{with} \quad \ctg(x,\check c_R) 
     := \min_{y \in \RR^m} \quad & \check c_R^\top y 
     + \tropindi{Ay+Bx\leq b}
    \end{aligned}
\end{equation}
where $\check p_R:= \bprob{\va c \in \relint(R)}$ and $\check c_R:=\besp{\va c \, | \, \va c \in \relint(R)}$ if $\check p_R>0$ and $\check c_R:=0$ if $\check p_R=0$.

In particular, if $\cR$ is a refinement of $\bigwedge_{\sigma \in \chcmplx{P}{\prjx}} -\cN_\sigma$,
\cref{eq:red_scen_det} holds for all $x \in \prj(P)$.
\end{theorem}
This is an exact quantization result, since~\eqref{eq:red_scen_det}
shows that $V(x)$ coincides with the value function of a
second stage problem with a cost distribution supported by the finite
set $\{\check c_R\mid R\in \cR\}$.

\begin{proof}
Let $\sigma \in \chcmplx{P}{\prj}$ and $x \in \relint(\sigma)$
 then, by definition, $\cN(P_x)=\cN_\sigma$.

For $R \in \cR$, there exists one and only one $N \in -\cN_\sigma$ such that $\relint(R) \subset \relint(N)$, that we denote $N(R)$. 
Indeed, as $\cR$ is a refinement of $-\cN_\sigma$, there exists at least one, and as $-\cN_\sigma$ is a polyhedral complex it is unique.

By \cref{lem:value_ctg}, under \cref{as:supp_vac} and since $x \in \prj(P)$,

\begin{subequations}
\begin{align*}
\ectg(x)&=\Besp{\sum_{N \in \cN(P_x)}\indi{c \in-\relint N} c^\top y_N(x)}
\\
&=  \Besp{\sum_{N \in -\cN_\sigma}\;\sum_{R \in \cR | \relint(R) \subset \relint(N)}\indi{\va c \in\relint R} \; \va c^\top y_N(x) }
&\text{ by the partition property}\\
&= \sum_{R \in \cR}\Ecind{\relint R} y_{N(R)}(x)
&\text{by linearity}\\
&= \sum_{R \in \cR} \check p_R \check c_R^\top y_{N(R)}(x) \\
& =\sum_{R \in \cR} \check p_R\min_{y \in \RR^m} \check c_R^\top y + \tropindi{Ay+Bx\leq b}
\end{align*}
\end{subequations}
the last equality is obtained by definition of $y_{N(R)}(x)$ as $\check c_R \in N(R)$, 
which leads to \cref{eq:red_scen_det}.

\qedb
 \end{proof}

Note that $\cR = \bigwedge_{\sigma \in {\chcmplxmax{P}{\prjx} }} -\cN_\sigma$ satisfies the condition of \cref{thm:quantization_2stage}
since if $\tau$ is a face of $\sigma$ in the chamber complex, $\cN_\sigma$ refines $\cN_\tau$
by \cite[Lemma 2.2]{rambau1996projections}.

\begin{cor}
	Under \cref{as:supp_vac},
let $x \in \prj(P)$ and $\sigma \in \chcmplx{P}{\prj}$ such that $x \in \relint(\sigma)$,
then for every refinement $\cR$ of $-\cN_\sigma$, the subgradient of $\ectg$ at point $x$
is given by the Minkowksi sum
\begin{equation}
\label{eq:subgradient_det}
      \partial \ectg(x)=  \sum_{R \in \cR} \check p_R BD(x,\check c_R)
\end{equation}
where $D(x,c):=\argmax \ba{(Bx-b)^\top\mu : A^\top \mu =- c, \mu\geq 0}$, 
 $\check p_R:= \bprob{\va c \in \relint(R)}$ and $\check c_R:=\besp{\va c \, | \, \va c \in \relint(R)}$ if $\check p_R>0$ and $\check c_R:=0$ if $\check p_R=0$.

In particular, if $\cR$ is a refinement of $\bigwedge_{\sigma \in \chcmplx{P}{\prjx}} -\cN_\sigma$, the subgradient formula \cref{eq:red_scen_det} holds for all $x \in \prj(P)$.
\label{cor:subgradient_2stage}
\end{cor}

\begin{proof}
	This is a consequence of the quantization result \cref{thm:quantization_2stage} and the formula of the subgradient of the expected cost-to-go function (see e.g.\ \cite[2.36]{shapiro2014lectures}), taking into account the form of the dual problem \cref{pb:dual_2nd-stage}.
\end{proof}

\begin{theorem}[Affine regions]
        For all distributions of $\va c$ satisfying~\cref{as:supp_vac},
the expected cost-to-go function $\ectg$ is affine on each cell of the chamber complex
        $\chcmplx{P}{\prjx}$.
	\label{thm:affine_chcmplx_2stage}
\end{theorem}

\begin{proof}

	Let $\sigma \in \chcmplx{P}{\prj}$. 
	We show that for every $c \in -\cone(A^\top)$, $x \mapsto \ctg(x,c)$ is an affine function on $\relint(\sigma)$.
By \cite{rambau1996projections} (see Lemma 2.1 (iii) and the comment after this lemma),
there exists a
unique minimal
face of $P$, denoted $F_{\sigma,c}$ that contains $\{x\}\times P_x^c$ for all $x \in \relint(\sigma)$.
Thus, for $x \in \relint(\sigma)$, $\ctg(x,c)=c^\top y$ for every $y \in P_x^c$ or equivalently for every $(x,y) \in F_{\sigma,c}$.
Since, $F_{\sigma,c}$ is a face, there exists an affine selection $x\mapsto y(x)$ on $\relint(\sigma) \subset \prj(F_{\sigma,c})$ such that $\bp{x,y(x)} \in F_{\sigma,c}$.
Then, $\ctg(\cdot,c)=c^\top y(\cdot)$ is affine on $\relint(\sigma)$.
By the quantization result, $\ectg(\cdot)=\sum_{R \in \cR} \check p_R \ctg(\cdot,\check c_R)$ is affine on $\relint(\sigma)$.

\end{proof}

\begin{remark}
It follows from this theorem that, for all $x\in \prjx(P)$,
	\begin{equation}
	\ectg(x) = \! \max_{\sigma \in {\chcmplxmax{P}{\prjx}}} \! \alpha_\sigma^\top x + \beta_\sigma  
	\text{ with }
	\alpha_\sigma = \! \sum_{N \in -\cN_\sigma}\! B^\top \mu_{\sigma}(\check c_N) \text{ and }
		\beta_\sigma = \! \sum_{N \in -\cN_\sigma}\! -b^\top \mu_{\sigma}(\check c_N)
\label{eq:formula_ectg_alpha_beta}
\end{equation}

where $\mu_{\sigma}(\check c_N) \in D(x,\check c_N)$ (defined in \cref{cor:subgradient_2stage}) for $x \in \relint(\sigma)$.
\end{remark}

%% file: sections/basis_decomposition.tex
The exact $2$-stage quantization theorem, \Cref{thm:quantization_2stage},
provides a polyhedral description of the expected cost-to-go function without recourse.
We next present
a combinatorial interpretation, through regular subdivisions and triangulations,
based on a result of Walkup and Wets, describing the piecewise
linear behavior of the value function of a deterministic
linear program. This provides further insight, and
also leads to an alternative way to prove~\Cref{thm:quantization_2stage}. 
However, this section is not used in the rest of the paper.

We first recall basic notions,
concerning the secondary fan, regular subdivisions and triangulations,
referring to the monograph of
De Loera, Rambau and Santos \cite{de2010triangulations}
for background. 

Let us denote by $(a_i)_{i\in [q]}$ the rows of the matrix $A$,
and let us choose $b\in \RR^q$. 
We shall think of $(a_i^\top)_{i\in [q]}$ as a {\em vector configuration}
in $\RR^m$, and $b$ as a {\em height} vector:
for each $i\in [q]$, we draw the point $(a_i^\top,b_i)\in \RR^m\times \RR$.
We now consider the convex hull $E$ of the points $(a_i^\top,b_i)$
in $\RR^m\times \RR$. The {\em geometric regular subdivision}
induced by the height vector $b$ is the polyhedral complex
defined as the projection onto $\RR^m$ of the lower faces
of the polyhedron $E$ (i.e., the faces of $E$ that a have a normal vector with a negative
ultimate coordinate). This geometric notion can be translated
in terms of a notion of {\em combinatorial} regular subdivision,
representing a face by the set of points that it contains.
These notions are formalized by the following definition.

\begin{definition}[Regular subdivisions, triangulations and secondary fan]
	Let us denote by $(a_i)_{i\in [q]}$ the rows of the matrix $A$,
        and let $b\in \RR^q$.
	The (combinatorial) \emph{regular subdivision} of the configuration of vectors $(a_i^\top)_{i\in [q]}$ induced by the height vector $b$ is the collection $\cS(A^\top,b)$ of
        subsets of $[q]$ such that
	\begin{equation}
		\cS(A^\top,b):= \ba{I \subset [q] \, |\, \exists y \in \RR^m,  a_i y=  b_i, \; \forall i \in I \text{ and } a_j y < b_j , \; \forall j \notin I} \enspace .
	\end{equation}
	A regular subdivision is a \emph{regular triangulation} when every
        set $I\in \cS(A^\top,b)$ yields an independent family of vectors $(a_i^\top)_{i \in I}$.

	When $A$ is fixed, the equivalence classes of $b\sim b' \iff \cS(A^\top,b)=\cS(A^\top,b')$ are relatively open cones. The collection of the closures of these cones constitutes a finite polyhedral fan, called the \emph{secondary fan}, and denoted by $\secfan{A^\top}$.	
\end{definition}

The geometry of the cost-to-go function $\ctg(x,c)$, for a deterministic $c$,
can be understood through the {\em basis decomposition theorem}
of Walkup and Wets \cite{walkup1969lifting}, see also the paper by 
Sturmfels and Thomas  \cite{sturmfels1997variation} for a more recent discussion.
In particular, Lemma 1.4 in  \cite{sturmfels1997variation}, applied to the dual problem \cref{pb:dual_2nd-stage} gives the following result.

\begin{lemma}[Basis decomposition and subdivision, see \cite{sturmfels1997variation}]
The set $D(x,c)$ of optimal solutions of the second stage dual problem \cref{pb:dual_2nd-stage},
defined in~\cref{cor:subgradient_2stage}, satisfies 
\begin{equation}
D(x,c)= \{ \mu \in \RR^q \, |\, A^\top \mu =-c, \mu\geq 0, \exists I \in \cS(A^\top,b-Bx),\, \supp(x) \subset I \} \enspace .
\end{equation}
	In particular, if $x$ and $x'$ belong to $\prj(P)$, the two following assertions are equivalent:
		\begin{itemize}
			\item[(i)] $b-Bx$ and $b-Bx'$ lie in the same relatively open cone of the secondary fan $\secfan{A^\top}$.
			\item[(ii)] For every $c \in -\cone(A^\top)$, $D(x,c)=D(x',c)$.
		\end{itemize}
	Futhermore, for $x \in \prj(P)$, the following assertions are equivalent
		\begin{itemize}
			\item[(iii)] $b-Bx$ lies in the interior of a maximal cone of the secondary fan $\secfan{A^\top}$.
			\item[(iv)] For every $c \in -\cone(A^\top)$, $D(x,c)$ is a singleton.
			\item[(v)] $\cS(A^\top,b-Bx)$ is a triangulation.
		\end{itemize}
	\label{lem:basis_decomposition}
\end{lemma}

Moreover, if we denote by $a$ the affine function such that $a(x):=b-Bx$, we can show that
\begin{align}
\chcmplx{P}{\prj}=a^{-1}\bp{\secfan{A^\top}}
\enspace .
\end{align}

Together with this fact, the basis decomposition theorem may allow us to retrieve our previous results.
Nevertheless, the proof of the equivalence between the chamber complex and the affine section of the secondary fan appearing above is technical.
We believe that the proof presented in \cref{ssec:finite_scenario} enlightens better the geometric and polyhedral structure of the expected cost-to-go function.

Recall that $\bigwedge_{\sigma \in \chcmplx{P}{\prjx}} \cN_\sigma$ appears in the exact quantization result for all $x$ in \cref{thm:quantization_2stage}.
This fan equals the chamber complex $\chcmplx{\cN(P)}{\xyprjy}$ of the normal fan $\cN(P)$ of the coupling constraint polyhedron along the projection $\xyprjy: (x,y)\mapsto y$.
Moreover, it is also the normal fan of the {\em fiber polyhedron}
$\Sigma\bp{P,\prjx(P)}$ defined in \cite{billera1992fiber}.
This is no coincidence, as the dual formulation of the 2-stage problem can be understood thanks to a (simple) generalization of the fiber polytope of \cite{billera1992fiber}.
However, to extend this interpretation to the multi-stage setting,
we need a more subtantial generalization of 
fiber polytopes, taking into account nonanticipativity constraints
and the nested structure of the control problem. We discuss
such a generalization in a subsequent work. In the next section, we develop
a direct approach to the multistage problem, in terms of chamber
complexes.

%% file: sections/extensions.tex
In this section, we show that the exact quantization result
established above for a general cost distribution and deterministic constraints 
carries over to the case of stochastic constraints with finite support and then to multistage programming.

We denote by $\xyprjx$ for the projection from 
$\RR^n \times \RR^m$ to $\RR^n$ defined by $\xyprjx(x',y')=x'$.
The projections $\xyzprjxy$, $\xyzprjx$, $\yzprjy$, $\prj^{x_{t-1},z}_{x_{t-1}}$ are defined accordingly. 
Note that in the notation $\xyzprjx$, $x$, $y$ and $z$ are part of the notation and not parameters.

\subsection{Propagating chamber complexes through Dynamic Programming}
\label{ssec:ext_sto_cons}
	
	We next show that chamber complexes are propagated through dynamic programming in a way that is \emph{uniform} with respect to the cost distribution.
	This is a key tool to extend the exact quantization theorem to the multistage setting.
Note that the proof of \cref{thm:quantization_2stage} cannot be extended to the multistage setting as, in this case,
the extensive form requires non-anticipativity constraints that cannot be tackled directly.

Recall that, for a polyhedron $P$ and a vector $\psi$,
we denote $P^\psi:= \argmin_{x \in P} \psi^\top x$.
Let $f$ be a polyhedral function on $\RR^d$, with a slight abuse of notation we denote 
$\epi(f)^{\psi,1} = \argmin_{(x,z) \in \epi(f)} \psi^\top x + z $.
We denote $\cFlow\bp{\epi(f)}:=\{\epi(f)^{\psi,1} \;|\; \psi \in \RR^d \}$ the set of \emph{lower faces}  of $\epi(f)$. 
The collection of projections (on $\RR^d$) of lower faces of $\epi(f)$ is the coarsest polyhedral complex such that $f$ is affine on each of its cells (see \cite[Chapter 2]{de2010triangulations}). 
Moreover, we have
\begin{equation}
\label{eq:equivalence_proj_E_psi_argmin}
    \prj_{\RR^d}\bp{(\epi(f)^{\psi,1}} = \argmin_{x \in \RR^d} \psi^\top x + f(x)
\end{equation}

\begin{lemma}
\label{lem:DP_chamber_complex}
    Let $\rctg$ be a polyhedral function on $\RR^m$ and $\rafreg :=\yzprjy\Bp{\cFlow\bp{\epi(\rctg)}}$ a coarsest
    polyhedral complex such that $\rctg$ is affine on each element of $\rafreg$. 
Let $\xi = (A,B,b)$ be fixed
and \cref{as:supp_vac} holds.
Define, for all $x \in \RR^n$
\begin{subequations}
\label{eq:ctg_with_recourse}
    \begin{align}
    Q(x,y) & :=   \rctg(y) + \tropindi{Ay + Bx \leq b} \\
\ectg(x)& := \besp{\min_{y \in \RR^m} \va c^\top y + Q(x,y) }  
\end{align}
\end{subequations}

Let $\cV:=\chcmplx{\cF(P)\wedge(\RR^n\times \rafreg)}{\xyprjx} \subset 2^{\RR^{n}}$
 with $P:=\{(x,y)\;|\;Ay+Bx\leq b \}$.

Then, $\cV\refleq \chcmplx{\epiQ}{\xyzprjx}$ and $\ectg$ is a polyhedral function which is affine on each element of $\cV$.
\end{lemma}

\input{tikz/dyn_prog_chamber_complex}

\begin{proof}

We have $\epiQ=\bp{\RR^n\times \epi(\rctg)}\cap (P\times \RR)\subset \RR^{n+m+1}$.
Since 
$$\ectg(x)= \besp{\min_{y \in \RR^m,z\in\RR} \va c^\top y + z + \tropindi{(x,y,z)\in \epiQ} } ,$$
by \cref{thm:affine_chcmplx_2stage} applied to the problem with variables $(y,z)$ and the coupling polyhedron $\epiQ$,
$\ectg$ is  a polyhedral function affine on each element of $\chcmplx{\epiQ}{\xyzprjx}$. 
We now show that $\cV\refleq \chcmplx{\epiQ}{\xyzprjx}$.
As $\epiQ$ is the epigraph of a polyhedral function, 
$\afregQ:=\xyzprjxy\bp{\cFlow(\epiQ)}%
\subset 2^{\RR^{n+m}}$ is a polyhedral complex.

Let $x_0 \in \xyzprjx(\epiQ)$, using notation of \cref{def:chmbrcplx},
\begin{align}
\sigma_{\epiQ,\xyzprjx}(x_0)
& :=\bigcap_{F\in \cF(\epiQ) \,\st\, x_0 \in \xyzprjx (F)} \xyzprjx(F) \\
& =\bigcap_{F\in \cFlow(\epiQ) \,\st\, x_0\in \xyzprjx(F)} \xyzprjx(F) \\
& =\bigcap_{F'\in \afregQ \,\st\, x_0\in \xyprjx(F')} \xyprjx(F')
=:\sigma_{\afregQ,\xyprjx}(x_0)
\end{align}

Indeed, as $\epiQ$ is an epigraph of a polyhedral function, if 
$F \in \cF(\epiQ)$ such that $x_0 \in \xyzprjx(F)$ 
then 
there exists 
$G \in \cFlow(\epiQ)$ such that $G \faceleq F$ and 
$x_0 \in \xyzprjx(G)$, allowing us to go from the first to second equality.
The third equality is obtained by setting $F'=\xyzprjxy(F)$.
Thus, $\chcmplx{\epiQ}{\xyzprjx}=\chcmplx{\afregQ}{\xyprjx}$.

We now show that $\cF(P)\wedge(\RR^n\times \rafreg) \refleq \afregQ$.
Let $G \in \cF(P)\wedge \bp{\RR^n \times \rafreg}$. There exist $\sigma \in \rafreg$ and $F \in \cF(P)$ such that $G=F\cap(\RR^n\times \sigma)$. %
By definition of $\cFlow$ , there exists  $\psi \in \RR^m$ such that $\sigma=\yzprjy\bp{\epi(\rctg)^{\psi,1}}$. 
We show that $G \subset \xyzprjxy(\epiQ^{0,\psi,1}) \in \afregQ$.
Indeed, let $(x,y) \in G=F\cap\bp{\RR^n \times \yzprjy (\epi(\rctg)^{\psi,1})}$. 
 We have $(x,y) \in F \subset P$ such that $y \in \argmin_{y' \in \RR^m} \ba{   \psi^\top y'+\rctg(y') }$. 
Which implies that $(x,y) \in \argmin \ba{   \psi^\top y'+ \rctg(y') \; | \; (x^\prime,y^\prime)\in P}$. 
This also reads, by \cref{eq:equivalence_proj_E_psi_argmin}, as $(x,y) \in \xyzprjxy(\epiQ^{0,\psi,1})$. 
Thus, $G \subset \xyzprjxy(\epiQ^{0,\psi,1}) \in \afregQ$ leading to $\cF(P)\wedge(\RR^n\times \rafreg) \refleq \afregQ$.
Finally, by monotonicity, \cref{lem:chmbrcplx_refinement} ends the proof. \qedb
\end{proof}

\begin{remark}
In \cref{lem:DP_chamber_complex}, the complex $\cV$ is independent of the distribution of $\va c$.
However, for special choices of $\va c$, $\ectg$ might be affine on each cell of a coarser complex than $\cV$.
For instance, if $\rctg=0$ and $\va c\equiv 0$, we have that 
$\ectg= \tropindi{\xyprjx(P)}$, $\ectg$ is affine on $\xyprjx(P)$.
Nevertheless, $\cV=\chcmplx{P}{\xyprjx}$ is generally finer than  $\cF\bp{\xyprjx(P)}$.
\end{remark}
\subsection{Exact quantization of MSLP}
	\label{ssec:ext_multi}

We next
show that the multistage program with arbitrary cost distribution is equivalent to a multistage program with independent, finitely distributed, cost distributions.
 Further, for all step $t$, there exist affine regions, independent of the distributions of costs, where $\ectg_t$ is affine.
\cref{as:supp_vac} is naturally extended to the multistage setting as follows

	\begin{hypo}
The sequence $(\va c_t,\va \xi_t)_{2\leq t\leq\horizon}$ is independent.\footnote{The results can be adapted to non-independent $\va \xi_t$ as long as $\va c_t$ is independent of $(\va{c}_\tau)_{\tau < t}$ conditionally on $(\va{\xi}_{\tau\leq t})$.} 
Further, for each $t\in\{2,\cdots,\horizon\}$, $\va \xi_t = ( \va A_t,\va B_t, \va b_t)$ is finitely supported, and  $\va c_t \in \spacevaintt$ is integrable with $\va c_t \in -\cone(\va A_t^\top)$ almost surely.
    \label{as:supp_vac_multi}
\end{hypo}

Note that \cref{as:supp_vac_multi} does not require independence between $\va c_t$ and $\va \xi_t$.
For $t\in[\horizon]$, and $\xi= (A,B,b) \in \supp(\va{\xi}_t)$ we define the coupling polyhedron 
\begin{equation}
    P_t(\xi):=\{(x_{t-1},x_t) \in \RR^{n_{t-1}}\times\RR^{n_{t}} \; | \; Ax_t+ B x_{t-1} \leq b\},
\end{equation}
and consider, for $x_{t-1}\in\RR^{n_{t-1}}$,
\begin{equation}
\label{eq:def_Vtilde}
    \tilde \ectg_t(x_{t-1}|\xi):= \besp{\min_{x_t \in \RR^{n_t}} \va c_t^\top x_t +V_{t+1}(x_t)+\tropindi{A x_t+ B x_{t-1}  \leq b} \;|\; \va \xi_t=\xi}
.
\end{equation}
Then, the cost-to-go function $\ectg_t$ is obtained by
\begin{equation}
\label{eq:V_as_finite_sum_of_Vtilde}
    \ectg_t(x_{t-1}) = 
    \sum_{\xi \in \supp(\va \xi_t)} \bprob{\va \xi_t = \xi}
   \tilde \ectg_t(x_{t-1}|\xi)
\end{equation}

The next two theorems extend the quantization results of \cref{thm:quantization_2stage} to the multistage settings.

\begin{theorem}[Affine regions independent of the cost]
\label{thm:polyhedral_multistage}
Assume that $(\va \xi_t)_{t\in [\horizon]}$ is a sequence of independent, finitely supported, random variables.
We define by induction $\cP_{\horizon+1}:=\{\RR^{n_{\horizon}}\}$ and for $t \in \{2,\dots,\horizon \}$
\begin{subequations}
\begin{align}
\cP_{t,\xi}&:=\chcmplx{\RR^{n_t}\times \cP_{t+1} \wedge  \cF\bp{P_t(\xi)} }{\prj^{x_{t-1},x_t}_{x_{t-1}}}\\
\cP_t &:= \bigwedge_{\xi_t \in \supp \va \xi_t} \cP_{t,\xi}
\end{align}
\end{subequations}
 Then, for all costs distributions $(\va c_t)_{2\leq t\leq\horizon}$ 
such that $(\va c_t, \va \xi_t)_{2\leq t \leq\horizon}$ 
satisfies \cref{as:supp_vac_multi} and all $t\in  \{2,\dots,\horizon \}$, 
we have $\supp(\cP_t) = \dom(\ectg_t)$, and $\ectg_t$ is polyhedral and affine on each cell of $\cP_t$.
 \end{theorem}

\begin{proof}
We set for all $t \in  \{2,\dots,\horizon+1 \}$, $\cV_t:=\prj^{x_{t-1},z}_{x_{t-1}}\bp{\cFlow\bp{\epi(\ectg_{t})}}$ the affine regions of $\ectg_t$.
As $\ectg_{\horizon+1}\equiv 0$
 is polyhedral and affine on $\RR^{n_{\horizon}}$, we have $\cP_{\horizon+1}
 =\cV_{\horizon+1}$.
Assume now that for $t \in  \{2,\dots,\horizon \}$, $\ectg_{t+1}$ is polyhedral and  $\cP_{t+1}$ refines
 $\cV_{t+1}$ (i.e. $\ectg_{t+1}$ is affine on each cell $\sigma \in \cP_{t+1}$).

By \cref{lem:DP_chamber_complex}, $\tilde \ectg_t(\cdot|\xi)$, defined in \cref{eq:def_Vtilde},
is affine on each cell of 
$\chcmplx{\RR^{n_t}\times\cV_{t+1}\wedge \cF\bp{P_t(\xi)} }{\pi^{x_{t-1},x_t}_{x_{t-1}}}$
which is refined by 
$\cP_{t,\xi}=\chcmplx{\RR^{n_t}\times\cP_{t+1}\wedge \cF\bp{P_t(\xi)} }{\pi_{x_{t-1}}^{x_{t-1},x_t}}$ by induction hypothesis  and \cref{lem:chmbrcplx_refinement}.
Thus, by \cref{eq:V_as_finite_sum_of_Vtilde}, $\ectg_t$ is affine on each cell of $\cP_t$.
In particular, $\ectg_t$ is polyhedral and $\cP_t := \bigwedge_{\xi_t \in \supp \va \xi_t} \cP_{t,\xi}$ refines 
$\cV_t$.
Backward induction ends the proof. \qedb
\end{proof}

By \cref{lem:DP_chamber_complex}, we have that $\cP_{t,\xi} \refleq \chcmplx{\epi\bp{Q_t^\xi}}{\prj^{x_{t-1},x_t,z}_{x_{t-1}}}$ where $Q_{t}^\xi(x_{t-1},x_t):=V_{t+1}(x_t)+\tropindi{Ax_{t}+Bx_{t-1}\leq b_t}$.
In particular, consider $\sigma \in \cP_{t,\xi}$, then 
for all $x_{t-1} \in \relint(\sigma)$, all fibers $\epi(Q_t^\xi)_{x_{t-1}}$ are normally equivalent.
We can then define $\cN_{t,\xi,\sigma}:=\cN(\epi(Q_t^\xi)_{x_{t-1}})$ for an arbitrary $x_{t-1} \in \relint(\sigma)$.

The next result shows that we can replace the MSLP problem \cref{eq:def_ectg_multi} by an equivalent problem with a discrete cost distribution.

\begin{theorem}[Exact quantization of the cost distribution, Multistage case]
Assume that $(\va \xi_t)_{t\in [\horizon]}$ is a sequence of independent, finitely supported, random variables.
Then, for all costs distributions
such that $(\va c_t, \va \xi_t)_{2 \leq t \leq \horizon}$ 
satisfies \cref{as:supp_vac_multi},
for all $t\in [\horizon]$, all  $x_{t-1} \in \RR^{n_{t-1}}$ and all $\xi \in \supp(\va \xi_t)$,
we have a quantized version of \cref{eq:def_Vtilde}:
     \begin{equation}
    \tilde \ectg_t(x_{t-1}|\xi) =
    \sum_{N \in \cN_{t,\xi}} \check p_{t,N|\xi} \; \min_{x_t \in \RR^{n_t}} 
    \Ba{\check{c}_{t,N|\xi}^\top x_t  + \ectg_{t+1}(x_t)+
    \tropindi{ A x_{t} +  B x_{t-1}\leq b}}
    \end{equation}
    where $\cN_{t,\xi} := \bigwedge_{\sigma \in \cP_{t,\xi}} -\cN_{t,\xi,\sigma}$ and for all $\xi \in \supp(\va \xi_t)$ and $N \in \cN_{t,\xi}$ we denote
    \begin{subequations}
\label{eq:def_check_p_check_c}
    \begin{align}
        \check p_{t,N|\xi} := \;& \bprob{\va c_t \in \relint N \; | \; \va \xi_t = \xi}
        \\
        \check c_{t,N|\xi} := \;& 
        \begin{cases}
            \besp{\va c_t \; | \; \va c_t \in \relint N, \va \xi_t = \xi} & 
            \text{ if } \bprob{\va \xi_t = \xi, \va x \in \relint N} \neq 0 \\
            0 & \text{otherwise}
        \end{cases}
    \end{align}
    \end{subequations}
\label{thm:multistage_quantization}
\end{theorem}

\begin{proof}

Since $\tilde \ectg_t(x_{t-1}|\xi)= \besp{\min_{x_t \in \RR^{n_t},z\in\RR} \va c^\top x_t + z + \tropindi{(x_{t-1},x_t,z)\in \epi(Q_t^\xi)} } $ and $\cP_{t,\xi}$ refines $\chcmplx{\epi\bp{Q_t^\xi}}{\prj^{x_{t-1},x_t,z}_{x_{t-1}}}$, by applying \cref{thm:quantization_2stage} with variables $(x_t,z)$ and the coupling constraints polyhedron $\epi(Q_{t}^\xi)$,
we deduce that
the coefficients $(\check p_{t,N|\xi})_{N\in\cN_{t,\xi}}$ and $(\check c_{t,N|\xi})_{N\in\cN_{t,\xi}}$ satisfy
	\begin{equation}
\tilde \ectg_t(x_{t-1}|\xi)=\sum_{N \in \cN_{t,\xi}} \check p_{t,N|\xi}
\min_{x_t \in \RR^{n_t}, z \in \RR} \Ba{\check c_{t,N|\xi}^\top x_t +z 
+ \tropindi{(x_{t-1},x_t,z) \in \epi(Q^\xi_{t})}}
\end{equation}
as the deterministic coefficient before $z$ is equal to its conditional expectation.
\qedb
 \end{proof}

In particular, the MSLP problem is equivalent to a finitely supported MSLP as shown in the following result.

For $t_0 \in [\horizon -1]$, we construct
the scenario tree $\cT_{t_0}$ as follows.
A node of depth $t-t_0$ of $\cT_{t_0}$ is labelled by
a sequence
$(N_\tau,\xi_\tau)_{t_0 < \tau \leq t}$ where 
    $N_\tau \in \cN_{\tau,\xi_\tau}$ and $\xi_\tau \in \supp ( \va{\xi}_\tau)$.
In this way, a node of depth $t-t_0$ of $\cT_{t_0}$ keeps track of the sequence
    of realizations of the random variables $\va{\xi}_\tau$ for times $\tau$ between
    $t_0$ and $t$, and of a selection of cones in $\cN_{t,\xi_{t}}$ at the same times.
Note that, by the independence assumption, all the subtrees of $\cT_{t_0}$, starting from a node of depth $t-t_0$ are the same as $\cT_{t_0+t}$. We denote by $\operatorname{lv}(\cT_{t_0})$ the set of leaves of $\cT_{t_0}$.

\begin{cor}[Equivalent finite tree problem]
        Define the quantized probability cost $c_\nu:= \check c_{t,N_t |\xi_t}$ and probability $p_\nu:= \prod_{t_0 < \tau \leq t} p_{\xi_\tau} \check p_{\tau,N_\tau|\xi_\tau}$, for all nodes $\nu=(N_\tau,\xi_\tau)_{t_0 <\tau \leq t}$.
Then, the cost-to-go functions associated with~\cref{eq:def_mslp}
		are given by
	\begin{subequations}
	\begin{align}
		V_{t_0}(x_0) = \min_{(x_\nu)_{\nu \in \cT_{t_0}}}&  \quad \sum_{\nu \in \cT_{t_0}} p_{\nu} c_{\nu}^\top x_{\nu} \\
		 \st  
		 & \quad A x_\mu+  B  x_{\nu} \leq  b  &
                 \forall
                 \nu \in \cT_{t_0}\backslash \operatorname{lv}(\cT_{t_0}), \forall \mu \succeq \nu\enspace,
	\end{align}
        for all $2\leq t_0\leq T-1$.
        Here, $x_0$ is the value of $x$ at the root node of $\cT_{t_0}$,
        and the notation $\forall \mu = (\nu,N,A,B,b) \succeq \nu$ indicates that $\mu$ ranges over the set of children of $\nu$.

	\label{pb:quantized_2stage_extensive}
	\end{subequations}
\label{cor:quantization_scenario_tree}
\end{cor}

%% file: tikz/dyn_prog_chamber_complex.tex
\def\raycone{0.5}
\def\lengthcone{0.55}
\def\colorcone{green!40!white}
\def\colorbigepi{blue!80!black}
\def\colorepiQ{yellow!40!white}
\def\colorR{yellow!90!orange}
\def\coloraffQ{green!70!black}
\def\colorP{orange!70!white}
\def\colorcE{brown!85!black}
\def\colorvrtx{brown!40!black}
\def\coloraffV{violet}
\def\colorcomref{red}
\def\colorbigpoly{blue}
\def\colorcuttingplane{yellow!60!white}
\def\colordash{black!60!white}

\def\origin{0,0,0}
\def\scalesub{0.9}
\def\originy{-2.7}

\def\zepi{2.1}

\def\pxa{0.7} 	\def\pya{-0.1}
\def\pxb{0.4} 	\def\pyb{-1.4}
\def\pxc{1.2} 	\def\pyc{-2}
\def\pxd{1.7} 	\def\pyd{-0.8}
\def\pxe{1.3} 	\def\pye{-0.25}

\def\qyf{-2.2} \def\qzf{1.7}
\def\qyg{-1.8} \def\qzg{1.2}
\def\qyh{-1} \def\qzh{0.6}
\def\qyi{-0.4} \def\qzi{0.8}

\tdplotsetmaincoords{65}{12}

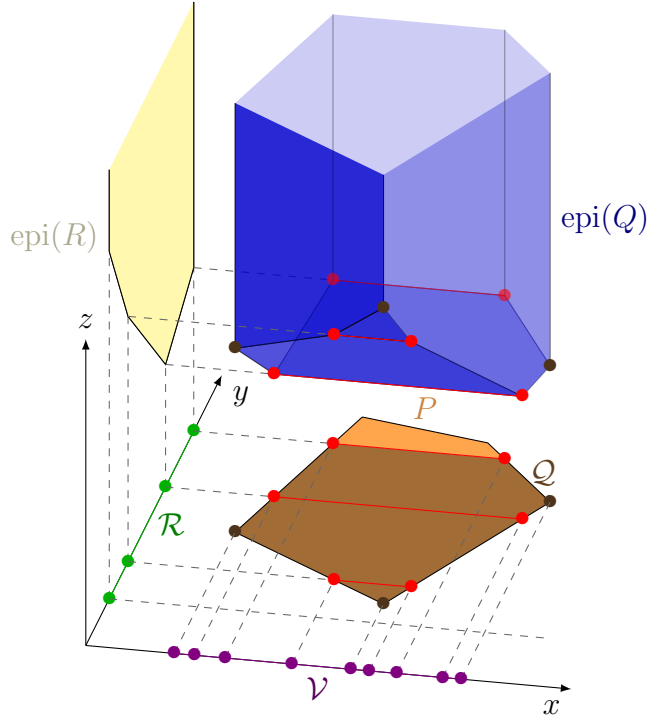
\begin{figure}[ht]
\centering
\begin{tikzpicture}[tdplot_main_coords,scale=\scaledynprogchcmplx]

\coordinate (a) at (\pxa,\pya,0);
\coordinate (b) at (\pxb,\pyb,0);
\coordinate (c) at (\pxc,\pyc,0);
\coordinate (d) at (\pxd,\pyd,0);
\coordinate (e) at (\pxe,\pye,0);

\coordinate (f) at (0,\qyf,\qzf);
\coordinate (g) at (0,\qyg,\qzg);

\coordinate (h) at (0,\qyh,\qzh);
\coordinate (i) at (0,\qyi,\qzi);
\coordinate (j) at (-0.5,-0.5,-0.5);

\coordinate (k) at (0,\qyf,\zepi);
\coordinate (l) at (0,\qyi,\zepi);

\pgfmathsetmacro\alphagh{(\qzg-\qzh)/(\qyg-\qyh)};
\pgfmathsetmacro\betagh{\qzg-\qyg*\alphagh};

\pgfmathsetmacro\alphaih{(\qzi-\qzh)/(\qyi-\qyh)};
\pgfmathsetmacro\betaih{\qzi-\qyi*\alphaih};

\pgfmathsetmacro\alphafg{(\qzf-\qzg)/(\qyf-\qyg)};
\pgfmathsetmacro\betafg{\qzf-\qyf*\alphafg};

\pgfmathsetmacro\alphaab{(\pxa-\pxb)/(\pya-\pyb)};
\pgfmathsetmacro\betaab{\pxa-\alphaab*\pya};
\pgfmathsetmacro\alphabc{(\pxb-\pxc)/(\pyb-\pyc)};
\pgfmathsetmacro\betabc{\pxb-\alphabc*\pyb};
\pgfmathsetmacro\alphacd{(\pxc-\pxd)/(\pyc-\pyd)};
\pgfmathsetmacro\betacd{\pxc-\alphacd*\pyc};
\pgfmathsetmacro\alphade{(\pxd-\pxe)/(\pyd-\pye)};
\pgfmathsetmacro\betade{\pxd-\alphade*\pyd};

\pgfmathsetmacro\abix{\alphaab*\qyi+\betaab};
\pgfmathsetmacro\abhx{\alphaab*\qyh+\betaab};
\pgfmathsetmacro\bcgx{\alphabc*\qyg+\betabc};
\pgfmathsetmacro\deix{\alphade*\qyi+\betade};
\pgfmathsetmacro\cdgx{\alphacd*\qyg+\betacd};
\pgfmathsetmacro\cdhx{\alphacd*\qyh+\betacd};

\draw[->] (0,\originy,0) -- (2.2,\originy,0) node[anchor=north east]{$x$};
\draw[->] (0,\originy,0) -- (0,0.2,0) node[anchor=north west]{$y$};
\draw[->] (0,\originy,0) -- (0,\originy,1.5) node[anchor=south]{$z$};

\draw[dashed,\colordash] (0,\qyf,0) -- (0,\qyf,\qzf);
\draw[dashed,\colordash] (0,\qyg,0) -- (0,\qyg,\qzg);
\draw[dashed,\colordash] (0,\qyh,0) -- (0,\qyh,\qzh);
\draw[dashed,\colordash] (0,\qyi,0) -- (0,\qyi,\qzi);

\fill[fill=\colorcE] (b) -- (c) -- (d)-- (\deix,\qyi,0) -- (\abix,\qyi,0);
\fill[fill=\colorP] (a) -- (e) -- (\deix,\qyi,0) -- (\abix,\qyi,0);
\draw (a) -- (b) -- (c) -- (d) -- (e) -- (a);
\draw[\colorvrtx] (b)node{$\bullet$};
\draw[\colorvrtx] (c)node{$\bullet$};
\draw[\colorvrtx] (d)node{$\bullet$};
\draw[\colorP!80!black] ($(a)!0.5!(e)$)node[anchor=south]{$P$};
\draw[\colorcE!60!black] ($(d)!0.5!(e)$)node[anchor=west]{$\afregQ$};

\draw (k) -- (f) -- (g) -- (h) -- (i) -- (l);
\filldraw[fill=\colorepiQ] (k) -- (f) -- (g) -- (h) -- (i) -- (l);
\draw[\colorepiQ!60!black] ($(k)!0.5!(f)$)node[anchor=north east]{$\epi (\rctg)$};
\draw[\coloraffQ] (0,\qyf,0)node{$\bullet$} -- (0,\qyg,0)node{$\bullet$} -- (0,\qyh,0)node{$\bullet$} -- (0,\qyi,0)node{$\bullet$};
\draw[\coloraffQ!70!black] ($(0,\qyg,0)!0.5!(0,\qyh,0)$)node[anchor=west]{$\rafreg$};

\draw[dashed,\colordash] (0,\qyf,0) -- (2,\qyf,0);
\draw[dashed,\colordash] (0,\qyg,0) -- (\bcgx,\qyg,0);
\draw[dashed,\colordash] (0,\qyh,0) -- (\abhx,\qyh,0);
\draw[dashed,\colordash] (0,\qyi,0) -- (\abix,\qyi,0);

\draw[dashed,\colordash] (0,\qyh,\qzh) -- (\abhx,\qyh,\qzh);
\draw[dashed,\colordash] (0,\qyi,\qzi) -- (\abix,\qyi,\qzi);

\draw[dashed,\colordash] (\pxb,\originy,0) -- (\pxb,\pyb,0);
\draw[dashed,\colordash] (\pxc,\originy,0) -- (\pxc,\pyc,0);
\draw[dashed,\colordash] (\pxd,\originy,0) -- (\pxd,\pyd,0);
\draw[dashed,\colordash] (\abhx,\originy,0) -- (\abhx,\qyh,0);
\draw[dashed,\colordash] (\abix,\originy,0) -- (\abix,\qyi,0);
\draw[dashed,\colordash] (\bcgx,\originy,0) -- (\bcgx,\qyg,0);
\draw[dashed,\colordash] (\cdgx,\originy,0) -- (\cdgx,\qyg,0);
\draw[dashed,\colordash] (\deix,\originy,0) -- (\deix,\qyi,0);
\draw[dashed,\colordash] (\cdhx,\originy,0) -- (\cdhx,\qyh,0);

\draw[\colorcomref](\bcgx,\qyg,0)node{$\bullet$} -- (\cdgx,\qyg,0)node{$\bullet$};
\draw[\colorcomref](\abhx,\qyh,0)node{$\bullet$} -- (\cdhx,\qyh,0)node{$\bullet$};
\draw[\colorcomref](\abix,\qyi,0)node{$\bullet$} -- (\deix,\qyi,0)node{$\bullet$};

\filldraw[\colorbigepi,draw=black,opacity=0.6]
 (\cdgx,\qyg,\qzg) -- (\bcgx,\qyg,\qzg) --(\pxb,\pyb,\alphagh*\pyb+\betagh) -- (\abhx,\qyh,\qzh) -- (\cdhx,\qyh,\qzh) -- (\cdgx,\qyg,\qzg);
 \filldraw[\colorbigepi,opacity=0.4]
 (\cdhx,\qyh,\qzh) -- (\abhx,\qyh,\qzh)node{$\bullet$} --  (\abix,\qyi,\qzi) -- (\deix,\qyi,\qzi) -- (\pxd,\pyd,\alphaih*\pyd+\betaih);
\filldraw[\colorbigepi,draw=black,opacity=0.7]
 (\cdgx,\qyg,\qzg) -- (\bcgx,\qyg,\qzg) --(\pxc,\pyc,\alphafg*\pyc+\betafg);

 \filldraw[\colorbigepi,draw=black,opacity=0.8]
 (\pxb,\pyb,\zepi)-- (\pxb,\pyb,\alphagh*\pyb+\betagh) --
 (\bcgx,\qyg,\qzg) --(\pxc,\pyc,\alphafg*\pyc+\betafg) -- (\pxc,\pyc,\zepi);

\filldraw[\colorbigepi,draw=black,opacity=0.3]
(\pxc,\pyc,\zepi)-- (\pxc,\pyc,\alphafg*\pyc+\betafg) --(\cdgx,\qyg,\qzg) -- (\cdhx,\qyh,\qzh) -- (\pxd,\pyd,\alphaih*\pyd+\betaih) -- (\pxd,\pyd,\zepi);

 \filldraw[\colorbigepi,draw=black,opacity=0.2]
(\pxd,\pyd,\zepi) -- (\pxd,\pyd,\alphaih*\pyd+\betaih) --(\deix,\qyi,\qzi) --(\deix,\qyi,\zepi);

 \filldraw[\colorbigepi,draw=black,opacity=0.2]
(\abix,\qyi,\zepi) -- (\abix,\qyi,\qzi) -- (\deix,\qyi,\qzi) --(\deix,\qyi,\zepi);
 \filldraw[\colorbigepi,draw=black,opacity=0.2]
(\abix,\qyi,\zepi) -- (\abix,\qyi,\qzi) -- (\abhx,\qyh,\qzh)--(\pxb,\pyb,\alphagh*\pyb+\betagh) -- (\pxb,\pyb,\zepi);

\draw[\colorbigepi!60!black] ($(\pxd,\pyd,\zepi)!0.5!(\pxd,\pyd,\alphaih*\pyd+\betaih)$)node[anchor=west]{$\epiQ$};

\draw[dashed,\colordash] (0,\qyg,\qzg) -- (\bcgx,\qyg,\qzg);

\draw[\colorcomref](\bcgx,\qyg,\qzg)node{$\bullet$} -- (\cdgx,\qyg,\qzg)node{$\bullet$};
\draw[\colorcomref](\abhx,\qyh,\qzh)node{$\bullet$} -- (\cdhx,\qyh,\qzh)node{$\bullet$};
\draw[\colorcomref,opacity=0.6](\abix,\qyi,\qzi)node{$\bullet$} -- (\deix,\qyi,\qzi)node{$\bullet$};

\draw[\colorvrtx] (\pxb,\pyb,\alphagh*\pyb+\betagh)node{$\bullet$};
\draw[\colorvrtx] (\pxc,\pyc,\alphafg*\pyc+\betafg)node{$\bullet$};
\draw[\colorvrtx] (\pxd,\pyd,\alphaih*\pyd+\betaih)node{$\bullet$};

\draw[\coloraffV] (\pxb,\originy,0)node{$\bullet$} -- (\abhx,\originy,0)node{$\bullet$} -- (\abix,\originy,0)node{$\bullet$} --  (\bcgx,\originy,0)node{$\bullet$} --
(\pxc,\originy,0)node{$\bullet$} -- (\cdgx,\originy,0)node{$\bullet$} -- (\deix,\originy,0)node{$\bullet$} -- (\cdhx,\originy,0)node{$\bullet$}  --
(\pxd,\originy,0)node{$\bullet$};
\draw[\coloraffV!70!black] ($(\pxb,\originy,0)!0.5!(\pxd,\originy,0)$)node[anchor=north]{$\cV$};

\end{tikzpicture}
\caption{An illustration of the proof of \cref{lem:DP_chamber_complex} : the epigraph $\epiQ$ of the coupling function in blue in the $(x,y,z)$ space, 
the epigraph of $\rctg$ in yellow in the $(y,z)$ plane, the affine regions $\rafreg$ of $\rctg$ in green on the $y$ axis, the coupling polyhedron $P$ in orange and brown in the $(x,y)$ plane, the polyhedral complex $\afregQ$ in red and brown in the $(x,y)$ plane and the chamber complex $\cV$ in violet on the $x$ axis.}
\end{figure}

%% file: sections/formulas.tex
	In this section, we show that, for three standard classes of distributions (uniform on a polytope, exponential, and Gaussian),
        the quantized costs $\check c_R$ and probabilities $\check p_R$
arising in the representation of the expected cost-to-go
function (\Cref{thm:quantization_2stage}),
        can be effectively computed.

The formulas are summed up in \cref{tab:synthesis_formulas_cost}.
They are detailed and established in
Sections~\ref{ss:unif}--\ref{ssec:gaussian}.
We provide these formulas for \emph{simplices} or 
\emph{simplicial cones} $S$ with $\dim(S)=\dim(\supp \va c)$.
This extends to any polyhedron $R$, through triangulation of $R \cap \supp(\va c)$ into simplices and simplicial cones $(S_k)_{k\in[l]}$. 
We then compute $\check p_R = \sum_{k=1}^l \check p_{S_k}$ and $\check c_R =  \sum_{k=1}^l \check p_{S_k} \check c_{S_k} / \check p_R$ if $\check p_R\neq 0$ and $\check c_R=0$ otherwise.

\begin{table}[ht]
\centering

{\footnotesize
\begin{tabular}{c|c|c|c}
&  Uniform %
& Exponential
& Gaussian  \\ \hline

$d\PP(c)$ 
& $\frac{\indi{c \in Q}}{\vol_{d}(Q)} d\cL_{\aff(Q)}(c)$
& $\frac{e^{\expvalpar^\top c} \indi{c \in K} }{\expval{K}}d\cL_{\aff(K)}c$
& $ \frac{e^{-\frac{1}{2}c^\top M^{-2}c}}{(2\pi)^\frac{m}{2}\det M}dc$
 \\\hline 

$\supp \va c$
& Polytope : $Q$
& Cone : $K$
& $\RR^m$
\\ \hline

$\check p_S$
& $\displaystyle \frac{\vol_{d}(S)}{\vol_{d} (Q)}$ 
&  $ \displaystyle \frac{|\det(\ray(S))|}{\expval{K}} \! \! \! \prod\limits_{r\in \ray(S)} \frac{1}{-r^\top \expvalpar}$
&  $\ang{M^{-1}S}$
\\ \hline

$\check c_S$
& $\frac{1}{d}\sum_{v \in \vrtx(S)} v$
& $\Big(\sum_{r\in \ray(S)} \frac{-r_i}{r^\top \expvalpar}\Big)_{i \in [m]}$
& $\frac{\sqrt 2 \Gamma(\frac{m+1}{2})}{\Gamma(\frac{m}{2})} M\sphcentr{S\cap \cansph{m-1}}$
\end{tabular}
}
\caption{Probabilities $\check p_S$ and expectations $\check c_S$ arising from different cost distributions over simplicial cones or simplices $S \subset \supp(\va c)$ with $\dim S = \dim (\supp \va c)$, where $\cL_A$ is the Lebesgue measure on an affine space $A$.}
\label{tab:synthesis_formulas_cost}
\end{table}

\subsection{Uniform distributions on polytopes}\label{ss:unif}

The \emph{volume} of a polytope $Q \subset \RR^m$ is the volume of $P$ seen as a subset of the smallest affine space $\aff(Q)$ it lives in.
The volume of a full dimensional simplex $S$ in $\RR^d$ with vertices $v_1,\dots,v_{d+1}$ is given by 
$\vol(S) = \frac{1}{n!}|\det(v_1-v_{d+1},\cdots, v_d- v_{d+1})|$,
see for example \cite{gritzmann1994complexityii} 3.1.
The \emph{centroid} of a non-empty polytope $Q \subset \RR^m$ is 
$\centr{Q}
:= \frac{1}{\vol Q} \int_Q y d\cL_{\aff Q}(y)$.
For instance, the centroid of a simplex $S$ of (non necessary full) dimension $d$ is the equibarycenter of its vertices :
$\centr{S}= \frac{1}{d+1} \sum_{v \in \vrtx (S)} v$.

Assume now that $Q$ is a polytope of dimension $\dimQ$,
and that $\va c$  is uniform on $Q$.
Let $S \subset Q$ be a simplex with $\dim(S)=\dim(Q)$, then we have
\begin{equation}
\label{eq:prob_uniform_poly}
	\check p_S= \frac{\vol_\dimQ{S}}{\vol_\dimQ{Q}}
	\qquad \text{and} \qquad
	\check c_S
 = \frac{1}{\dimQ+1} \sum_{v \in \vrtx(S)} v \enspace .
\end{equation}

\subsection{Exponential distributions on cones}
\label{ssec:exponential_distr}
Let $P$ be a polyhedron  and $\expvalpar \in  \relint \bp{(\rc P)^\circ}$.
We denote by $\expval{P} := \int_P e^{\expvalpar^\top c} d\cL_{\aff(P)}(c)$ the \emph{exponential valuation} of $P$ with parameter $\expvalpar$.

\begin{proposition}[Brion's formula \cite{brion1988points}]
Let $S$ be a full dimensional simplicial cone,
and let $\ray(S)$ be a square matrix
whose columns are obtained by selecting precisely one element
in every extreme ray of $S$, so that $S=\cone\bp{\ray(S)}$.
Then for any $\expvalpar \in \relint S^\circ$, the exponential valuation of $S$ is given by
\begin{equation}
\expval{S}= |\det(\ray(S))| \prod_{r \in \ray(S)} \frac{1}{-\expvalpar^\top r}
\enspace . \label{eq:brion_formula}
\end{equation}
\end{proposition}

Let $K$ be a (non necessarily simplicial) polyhedral cone and $\expvalpar \in \relint K^\circ$ a vector. 
Assume that $\va c$ has the following exponential density :
\begin{equation}
d\PP(c) := e^{\expvalpar^\top c} \indi{c \in K}  \frac{1}{\expval{K}}d\cL_{\aff(K)}(c)
\label{e-def-exptype}
\end{equation}
Let $S \subset K$ be a simplicial cone with $\dim S = \dim K$,
by Brion's formula \eqref{eq:brion_formula},

\begin{equation}
\check p_S
 = \frac{\expval{S}}{\expval{K}} = \frac{1}{\expval{K}} |\det(\ray(S))| \prod_{r\in \ray(S)} \frac{1}{-r^\top \expvalpar}
 \label{eq:prob_expval}
\end{equation}
Further, 
\begin{equation}
\check p_S \check c_S= 
\Ecind{S}
 = \frac{1}{\expval{K}} \int_{S} c e^{\expvalpar^\top c}dc
 = \frac{\nabla \Phi_{S} (\expvalpar)}{\expval{K}} \enspace .
\end{equation}
By computing explictly the latter gradient,
dividing by $\check{p}_S$, and simplifying,
we obtain:
\begin{equation}
\check c_S=\Big(\sum_{r \in \ray(S)} \frac{-r_i}{r^\top \expvalpar}\Big)_{i \in [m]}
\enspace .
\label{eq:esp_expval}
\end{equation}

\subsection{Gaussian distributions}
\label{ssec:gaussian}

The \emph{solid angle} 
of a pointed cone $K \subset \RR^d$ is defined as the normalized volume
of its intersection with the unit ball $\canball{d}$,
i.e.: $\ang{K}:={\vol_d (K \cap \canball{d})}/{\vol_d \canball{d}}$.
Recall that $\vol_d \canball{d}={\pi^\frac{d}{2}}/{\Gamma(\frac{d}{2}+1)}$ with $\Gamma$ the Euler gamma function, and that (\cite{ribando2006measuring})
for any function $f : \RR^m \to \RR$ invariant under rotations around the origin and any pointed cone $K \subset \RR^m$, we have 
$
\ang{K} \int_{\RR^m} f = \int_K f
$.

Let $\va c$ be a non-degenerate, centered, Gaussian random variable of variance $M^2$, where $M$ is a symmetric positive definite matrix.
Then, if $K$ is a polyhedral cone, we have
\begin{align*}
\check p_K
&=\int_{M^{-1}K} \frac{e^{-\frac{1}{2}\mynorm{c}_2^2}}{(2\pi)^\frac{m}{2}}dc
=\ang{M^{-1}K}
\end{align*}

We shall use the notion of {\em spherical centroid} $\sphcentr{U}$
for a measurable subset $U$ included in the unit sphere. It is defined
as the barycenter of the elements of $U$ with respect
to the uniform measure on the sphere. Note that the spherical centroid
does not belong to the sphere, unless $U$ is trivial. We have
\begin{align*}
\check p_K \check c_K
&=  \int_{M^{-1}K} Mc\frac{e^{-\frac{1}{2}\mynorm{c}_2^2}}{(2\pi)^\frac{m}{2}}dc
= M \int_{\RR^+} r^m \frac{e^{-\frac{r^2}{2}}}{(2\pi)^\frac{m}{2}} dr   \int_{M^{-1}K \cap \cansph{m-1}} \varphi d\varphi \\
&=M\frac{\Gamma(\frac{m+1}{2})}{\sqrt 2 \pi^{\frac{m}{2}}} \vol_{m-1}(\cansph{m-1}) \ang{M^{-1}K} \sphcentr{M^{-1}K \cap \cansph{m-1}} \\
&= M\frac{\sqrt 2 \Gamma(\frac{m+1}{2})}{\Gamma(\frac{m}{2})}  \ang{M^{-1}K} \sphcentr{M^{-1}K \cap \cansph{m-1}}
\end{align*}

Similarly, one can get explicit formul\ae\
when $\va{c}$ is distributed uniformly on an ellipsoid,
or on the surface of an ellipsoid, or more generally,
when the distribution of $\va{c}$ is invariant
under the action of an orthogonal group.
Then, the quantized costs and probabilities $\check{c}_S$ and
$\check{p}_S$ are still given by solid angles and spherical centroids,
in a way similar to~\Cref{tab:synthesis_formulas_cost}.

%% file: sections/example.tex
We consider the following second-stage problem, with $n=1$ and $m=2$ :

\begin{equation}
    \ectg(x) = \EE\left[
    \begin{aligned}\min_{y \in \RR^2} \quad & \va{c}^\top y\\
    \st \quad & \mynorm{y}_1 \leq 1, \quad  y_1 \leq x \text{ and } y_2 \leq x
    \end{aligned}\right] \enspace . 
\end{equation}
 We apply our results, to provide an explicit representation of $V$.

\input{sections/tikzfigure}

The coupling polyhedron is 
$P=\{(x,y) \in \RR^n\times \RR^m\, | \, \mynorm{y}_1 \leq 1,\, y_i \leq x \quad \forall i \in [m]\}$
presented in \cref{fig:coupling_polyhedron},
and its V-representation is the collection of vertices
$
(0,-1,0)$,
$(-\frac{1}{2},-\frac{1}{2},-\frac{1}{2})$, 
$(0,0,-1)$,
$(1,1,0)$,
$(\frac{1}{2},\frac{1}{2},\frac{1}{2})$,
$(1,0,1)$ and the ray $(1,0,0)
$.
By projecting the different faces, we see that its projection is the half-line,
$
\prj(P)=[-\frac{1}{2},+\infty[
$
and its chamber complex is
$
\chcmplx{P}{\prj}$ is the collection of cells composed of
$\{-\frac{1}{2}\}$,
$[-\frac{1}{2},0]$,
$\{0\}$,
$[0,\frac{1}{2}]$,
$\{\frac{1}{2}\}$,
$[\frac{1}{2},1]$,
$\{1\}$,
$[1,+\infty)
$ as presented in \cref{fig:coupling_polyhedron}.
As there are 4 different maximal chambers, there are 4 different classes of normally equivalent fibers as shown in \cref{fig:different_fibers}.

We evaluate $\check c_N$  and $\check p_N$ for $N \in -\cN_\sigma$ using
the formulas of~\Cref{tab:synthesis_formulas_cost}.
For example, when $\va c$ is uniform on the centered ball for the $\infty$-norm of radius $R$, \cref{fig:compfansupp_ball_infty} shows the regions of which
the areas and centroids need to be computed.

\begin{figure}[ht]
\input{tikz/compfansupp_ball_infty}
\label{fig:compfansupp_ball_infty}
\end{figure}

\begin{figure}[!ht]
\begin{tikzpicture}[scale=2.5]
\input{tikz/graph_Q_example}
\end{tikzpicture}
\centering 
\caption{Graph of the cost-to-go function $\ectg$ for different distribution of the cost $\va c$ with $R=\theta=\gamma=1$.}
\label{fig:graph_ectg}
\end{figure}

\input{table/affine_values_ectg}

%% file: sections/tikzfigure.tex
\def\raycone{0.5}
\def\lengthcone{0.55}
\def\colorcone{green!40!white}
\def\colorboundcone{green!70!blue}
\def\colorpoly{blue!40!white}
\def\colorbigpoly{blue}
\def\colorcuttingplane{yellow!60!white}
\def\origin{0,0,0}
\def\scalesub{\scalefibers}

\tdplotsetmaincoords{65}{12}

\begin{figure}[ht]

\centering

\begin{tikzpicture}[tdplot_main_coords,scale=\scalecoupP]
\input{tikz/big_poly_cuts}
\end{tikzpicture}
\caption{The coupling polyhedron $P$ in blue, different cuts and fibers $P_x$ vertical in yellow, and its chamber complex $\chcmplx{P}{\prjx}$ in red on the bottom.}
\label{fig:coupling_polyhedron}
\end{figure}

\begin{figure}[ht]
\begin{subfigure}{.27\textwidth}
  \begin{tikzpicture}[scale=\scalesub]
\input{tikz/fiber-025}
\end{tikzpicture}
\centering
\caption{$x=-0.25,\sigma=[-0.5,0]$}
\end{subfigure}%
\begin{subfigure}{.23\textwidth}
\begin{tikzpicture}[scale=\scalesub]
\input{tikz/fiber+025}
\end{tikzpicture}
\centering
\caption{$x=0.25, \sigma=[0,0.5]$}
\end{subfigure}
\begin{subfigure}{.24\textwidth}
\begin{tikzpicture}[scale=\scalesub]
\input{tikz/fiber+075}
\end{tikzpicture}
\centering
\caption{$x=0.75, \sigma=[0.5,1]$}
\end{subfigure}
\begin{subfigure}{0.24\textwidth}
  \begin{tikzpicture}[scale=\scalesub]
\input{tikz/fiber+1}
\end{tikzpicture}
\centering
\caption{$x\geq1, \sigma=[1,+\infty)$}
\end{subfigure}

\caption{Fibers $P_x$ in blue and their normal fan $\cN(P_x)=\cN_\sigma$ in green for different $x \in \RR$}

\label{fig:different_fibers}

\end{figure}

%% file: tikz/big_poly_cuts.tex
\coordinate (o) at (\origin);
\coordinate (a) at ($(o)+(2,1,0)$);
\coordinate (b) at ($(o)+(2,0,1)$);
\coordinate (c) at ($(o)+(2,-1,0)$);
\coordinate (d) at ($(o)+(2,0,-1)$);

\coordinate (e) at ($(o)+(1,1,0)$);
\coordinate (f) at ($(o)+(1,0,1)$);
\coordinate (g) at ($(o)+(0,-1,0)$);
\coordinate (h) at ($(o)+(0,0,-1)$);

\coordinate (i) at ($(o)+(0.5,0.5,0.5)$);
\coordinate (j) at ($(o)+(-0.5,-0.5,-0.5)$);

\draw (a) -- (b) -- (c) -- (d) -- (a);
\filldraw[fill=\colorbigpoly,opacity=0.9] (a) -- (b) -- (c) -- (d) -- (a);

\draw[\colorbigpoly]  (2.3,0,0) node[anchor=north,scale=1.2]{$P$};

\fill[\colorbigpoly,opacity=0.1] (a) -- (e) -- (h) --(d) -- (a);
\fill[\colorbigpoly,opacity=0.7] (b) -- (f) -- (g) --(c) -- (b);

\fill[\colorbigpoly,opacity=0.1] (i) -- (e) -- (h) --(j) -- (i);
\fill[\colorbigpoly,opacity=0.3] (i) -- (f) -- (g) --(j) -- (i);

\fill[\colorbigpoly,opacity=0.05] (a) -- (b) -- (f) --(i) -- (e)--(a);
\fill[\colorbigpoly,opacity=0.1] (c) -- (d) -- (h) --(j) -- (g) -- (c);

\draw[dashed] (a) -- (e);
\draw (b) -- (f);
\draw (c) -- (g);
\draw (d) -- (h);

\draw[dashed] (e) -- (h);
\draw (f) -- (g);

\draw[dashed] (i) -- (e);
\draw (i) -- (f);
\draw (j) -- (g);
\draw (j) -- (h);

\draw (i) -- (j);

\draw[->] (a) -- ($(a)+ \raycone*(1,0,0)$);
\draw[->] (b) -- ($(b)+ \raycone*(1,0,0)$);
\draw[->] (c) -- ($(c)+ \raycone*(1,0,0)$);
\draw[->] (d) -- ($(d)+ \raycone*(1,0,0)$);

\draw[dotted,->] (o) -- ($(o)+(3,0,0)$) node[anchor=north east]{$x$};
\draw[dotted,->] (o) -- ($(o)+(0,2,0)$) node[anchor=north west]{$y_1$};
\draw[dotted,->] (o) -- ($(o)+(0,0,1.7)$) node[anchor=south]{$y_2$};

\def\cutx{-0.25}
\fill[\colorcuttingplane,opacity=0.5] ($(o)+(\cutx,1,1)$)node[color=black,opacity=1,anchor=north,scale=0.7]{$x=\cutx$}  -- (\cutx,1,-1) -- (\cutx,-1,-1) --(\cutx,-1,1) -- cycle ;
\draw (\cutx,\cutx,\cutx) -- (\cutx,\cutx,-1-\cutx)  ;
\draw (\cutx,\cutx,-1-\cutx) -- (\cutx,-1-\cutx,\cutx);
\draw (\cutx,-1-\cutx,\cutx) -- (\cutx,\cutx,\cutx);

\def\cutx{0.25}
\fill[\colorcuttingplane,opacity=0.5] (\cutx,1,1)node[color=black,opacity=1,anchor=north,scale=0.7]{$x=\cutx$} -- (\cutx,1,-1) -- (\cutx,-1,-1) --(\cutx,-1,1) -- (\cutx,1,1);
\draw (\cutx,-1,0) -- (\cutx,0,-1)  -- (\cutx,\cutx,\cutx-1)-- (\cutx,\cutx,\cutx)--(\cutx,\cutx -1,\cutx) -- cycle;

\def\cutx{0.75}
\fill[\colorcuttingplane,opacity=0.5] (\cutx,1,1)node[color=black,opacity=1,anchor=north,scale=0.7]{$x=\cutx$} -- (\cutx,1,-1) -- (\cutx,-1,-1)  --(\cutx,-1,1) -- (\cutx,1,1);
\draw (\cutx,0,-1)-- (\cutx,-1,0) -- (\cutx,\cutx-1,\cutx)-- (\cutx,1-\cutx,\cutx)--(\cutx,\cutx,1-\cutx)--(\cutx,\cutx,\cutx-1) -- cycle;

\def\cutx{1.5}
\fill[\colorcuttingplane,opacity=0.5] (\cutx,1,1)node[color=black,opacity=1,anchor=north,scale=0.7]{$x=\cutx$} -- (\cutx,1,-1) -- (\cutx,-1,-1)  --(\cutx,-1,1) -- (\cutx,1,1);
\draw (\cutx,-1,0) -- (\cutx,0,-1)  -- (\cutx,1,0)-- (\cutx,0,1) -- cycle;

\def\ychcmplx{-1,-1}
\draw[red] (-0.5,\ychcmplx) node{$\bullet$} -- (0,\ychcmplx)node{$\bullet$} -- (0.5,\ychcmplx)node{$\bullet$} -- (1,\ychcmplx)node{$\bullet$} ;
\draw[red,->] (1,\ychcmplx) -- (2.8,\ychcmplx);
\draw (-0.5,\ychcmplx) node[anchor=north,scale=0.7]{$x=-0.5$};
\draw (0,\ychcmplx) node[anchor=north,scale=0.7]{$x=0$};
\draw (0.5,\ychcmplx) node[anchor=north,scale=0.7]{$x=0.5$};
\draw (1,\ychcmplx) node[anchor=north,scale=0.7]{$x=1$};
\draw[red] (1.7,\ychcmplx) node[anchor=north,scale=1.2]{$\chcmplx{P}{\prj}$};

%% file: tikz/fiber-025.tex
\coordinate (a) at (-0.25,-0.75);
\coordinate (b) at (-0.75,-0.25);
\coordinate (c) at (-0.25,-0.25);

\fill[white] ($(0,-1) + (315:\raycone)$) -- (0,-1) -- ($(0,-1)+ (225:\raycone)$)  arc (225:315:\raycone);

\filldraw[fill=\colorpoly, draw=black] (a) -- (b) -- (c) -- (a);

\fill[\colorcone] ($(c) + \raycone*(0,1)$) -- (c) -- ($(c) + \raycone*(1,0)$) arc (0:90:\raycone);
\draw[thick,->,\colorboundcone] (c) -- ($(c) + \lengthcone*(0,1)$);
\draw[thick,->,\colorboundcone] (c) -- ($(c) + \lengthcone*(1,0)$);

\draw[thick,->,\colorboundcone]  ($ (a)!0.5!(c) $) -- ($(a)!0.5!(c) + \lengthcone*(1,0)$);

\fill[\colorcone] ($(a) + (225:\raycone)$) -- (a) -- ($(a)+ \raycone*(1,0)$)  arc (0:-135:\raycone);
\draw[thick,->,\colorboundcone] (a) -- ($(a) + (225:\lengthcone)$);
\draw[thick,->,\colorboundcone] (a) -- ($(a) + \lengthcone*(1,0)$);

\draw[thick,->,\colorboundcone]  ($ (a)!0.5!(b) $) -- ($(a)!0.5!(b) + (225:\lengthcone)$);

\fill[\colorcone] ($(b) + (225:\raycone)$) -- (b) -- ($(b)+ \raycone*(0,1)$)  arc (90:225:\raycone);
\draw[thick,->,\colorboundcone] (b) -- ($(b) + (225:\lengthcone)$);
\draw[thick,->,\colorboundcone] (b) -- ($(b) + \lengthcone*(0,1)$);

\draw[thick,->,\colorboundcone]  ($ (c)!0.5!(b) $) -- ($(c)!0.5!(b) + (0,\lengthcone)$);

\draw[thin,dashed,->] (0,0) -- (1.25,0) node[anchor=west] {$y_1$};
\draw[thin,dashed,->] (0,0) -- (0,1.25) node[anchor=south] {$y_2$};
\draw[thin,dashed] (0,0) -- (-1.25,0);
\draw[thin,dashed] (0,0) -- (0,-1.25);

%% file: tikz/fiber+025.tex
\coordinate (a) at (0,-1);
\coordinate (b) at (-1,0);
\coordinate (c) at (-0.75,0.25);
\coordinate (d) at (0.25,0.25);
\coordinate (e) at (0.25,-0.75);

\filldraw[fill=\colorpoly, draw=black] (a) -- (b) -- (c) -- (d) -- (e) -- (a);

\fill[\colorcone] ($(a) + (315:\raycone)$) -- (a) -- ($(a)+ (225:\raycone)$)  arc (225:315:\raycone);
\draw[thick,->,\colorboundcone] (a) -- ($(a) + (315:\lengthcone)$);
\draw[thick,->,\colorboundcone] (a) -- ($(a) + (225:\lengthcone)$);

\draw[thick,->,\colorboundcone]  ($ (a)!0.5!(b) $) -- ($(a)!0.5!(b) + (225:\lengthcone)$);

\fill[\colorcone] ($(b) + (225:\raycone)$) -- (b) -- ($(b)+ (135:\raycone)$)  arc (135:225:\raycone);
\draw[thick,->,\colorboundcone] (b) -- ($(b) + (225:\lengthcone)$);
\draw[thick,->,\colorboundcone] (b) -- ($(b) + (135:\lengthcone)$);

\draw[thick,->,\colorboundcone]  ($ (c)!0.5!(b) $) -- ($(c)!0.5!(b) + (135:\lengthcone)$);

\fill[\colorcone] ($(c) + (135:\raycone)$) -- (c) -- ($(c)+ (90:\raycone)$)  arc (90:135:\raycone);
\draw[thick,->,\colorboundcone] (c) -- ($(c) + (90:\lengthcone)$);
\draw[thick,->,\colorboundcone] (c) -- ($(c) + (135:\lengthcone)$);

\draw[thick,->,\colorboundcone]  ($ (c)!0.5!(d) $) -- ($(c)!0.5!(d) + (90:\lengthcone)$);

\fill[\colorcone] ($(d) + \raycone*(0,1)$) -- (d) -- ($(d) + \raycone*(1,0)$) arc (0:90:\raycone);
\draw[thick,->,\colorboundcone] (d) -- ($(d) + \lengthcone*(0,1)$);
\draw[thick,->,\colorboundcone] (d) -- ($(d) + \lengthcone*(1,0)$);

\draw[thick,->,\colorboundcone]  ($ (d)!0.5!(e) $) -- ($(d)!0.5!(e) + \lengthcone*(1,0)$);

\fill[\colorcone] ($(e)+ \raycone*(1,0)$) -- (e) -- ($(e) + (315:\raycone)$)  arc (-45:0:\raycone);
\draw[thick,->,\colorboundcone] (e) -- ($(e) + (315:\lengthcone)$);
\draw[thick,->,\colorboundcone] (e) -- ($(e) + \lengthcone*(1,0)$);

\draw[thick,->,\colorboundcone]  ($ (e)!0.5!(a) $) -- ($(e)!0.5!(a) + (315:\lengthcone)$);

\draw[thin,dashed,->] (0,0) -- (1.25,0) node[anchor=west] {$y_1$};
\draw[thin,dashed,->] (0,0) -- (0,1.25) node[anchor=south] {$y_2$};
\draw[thin,dashed] (0,0) -- (-1.25,0);
\draw[thin,dashed] (0,0) -- (0,-1.25);

%% file: tikz/fiber+075.tex
\coordinate (a) at (0,-1);
\coordinate (b) at (-1,0);
\coordinate (c) at (-0.25,0.75);
\coordinate (d) at (0.25,0.75);
\coordinate (e) at (0.75,0.25);
\coordinate (f) at (0.75,-0.25);

\filldraw[fill=\colorpoly, draw=black] (a) -- (b) -- (c) -- (d) -- (e) -- (f) --(a);

\fill[\colorcone] ($(a) + (315:\raycone)$) -- (a) -- ($(a)+ (225:\raycone)$)  arc (225:315:\raycone);
\draw[thick,->,\colorboundcone] (a) -- ($(a) + (315:\lengthcone)$);
\draw[thick,->,\colorboundcone] (a) -- ($(a) + (225:\lengthcone)$);

\draw[thick,->,\colorboundcone]  ($ (a)!0.5!(b) $) -- ($(a)!0.5!(b) + (225:\lengthcone)$);

\fill[\colorcone] ($(b) + (225:\raycone)$) -- (b) -- ($(b)+ (135:\raycone)$)  arc (135:225:\raycone);
\draw[thick,->,\colorboundcone] (b) -- ($(b) + (225:\lengthcone)$);
\draw[thick,->,\colorboundcone] (b) -- ($(b) + (135:\lengthcone)$);

\draw[thick,->,\colorboundcone]  ($ (c)!0.5!(b) $) -- ($(c)!0.5!(b) + (135:\lengthcone)$);

\fill[\colorcone] ($(c) + (135:\raycone)$) -- (c) -- ($(c)+ (90:\raycone)$)  arc (90:135:\raycone);
\draw[thick,->,\colorboundcone] (c) -- ($(c) + (90:\lengthcone)$);
\draw[thick,->,\colorboundcone] (c) -- ($(c) + (135:\lengthcone)$);

\draw[thick,->,\colorboundcone]  ($ (c)!0.5!(d) $) -- ($(c)!0.5!(d) + (90:\lengthcone)$);

\fill[\colorcone] ($(d) + (90:\lengthcone)$) -- (d) -- ($(d) + (45:\lengthcone)$) arc (45:90:\raycone);
\draw[thick,->,\colorboundcone] (d) -- ($(d) + (90:\lengthcone)$);
\draw[thick,->,\colorboundcone] (d) -- ($(d) + (45:\lengthcone)$);

\draw[thick,->,\colorboundcone]  ($ (d)!0.5!(e) $) -- ($(d)!0.5!(e) + (45:\lengthcone)$);

\fill[\colorcone] ($(e) + (45:\raycone)$) -- (e) -- ($(e) + (0:\raycone)$)  arc (0:45:\raycone);
\draw[thick,->,\colorboundcone] (e) -- ($(e) + (45:\lengthcone)$);
\draw[thick,->,\colorboundcone] (e) -- ($(e) +  (0:\lengthcone)$);

\draw[thick,->,\colorboundcone]  ($ (e)!0.5!(f) $) -- ($(e)!0.5!(f) + \lengthcone*(1,0)$);

\fill[\colorcone] ($(f) + (0:\raycone)$) -- (f) -- ($(f) + (-45:\raycone)$)  arc (-45:0:\raycone);
\draw[thick,->,\colorboundcone] (f) -- ($(f) + (-45:\lengthcone)$);
\draw[thick,->,\colorboundcone] (f) -- ($(f) +  (0:\lengthcone)$);

\draw[thick,->,\colorboundcone]  ($ (f)!0.5!(a) $) -- ($(f)!0.5!(a) + (315:\lengthcone)$);

\draw[thin,dashed,->] (0,0) -- (1.25,0) node[anchor=west] {$y_1$};
\draw[thin,dashed,->] (0,0) -- (0,1.25) node[anchor=south] {$y_2$};
\draw[thin,dashed] (0,0) -- (-1.25,0);
\draw[thin,dashed] (0,0) -- (0,-1.25);

\draw[\colorboundcone]  ($ (e)!0.5!(f) $) -- ($(e)!0.5!(f) + \lengthcone*(1,0)$);

%% file: tikz/fiber+1.tex
\coordinate (a) at (0,-1);
\coordinate (b) at (-1,0);
\coordinate (c) at (0,1);
\coordinate (d) at (1,0);

\filldraw[fill=\colorpoly] (a) -- (b) -- (c) -- (d)--(a);

\fill[\colorcone] ($(a) + (315:\raycone)$) -- (a) -- ($(a)+ (225:\raycone)$)  arc (225:315:\raycone);
\draw[thick,->,\colorboundcone] (a) -- ($(a) + (315:\lengthcone)$);
\draw[thick,->,\colorboundcone] (a) -- ($(a) + (225:\lengthcone)$);

\draw[thick,->,\colorboundcone]  ($ (a)!0.5!(b) $) -- ($(a)!0.5!(b) + (225:\lengthcone)$);

\fill[\colorcone] ($(b) + (225:\raycone)$) -- (b) -- ($(b)+ (135:\raycone)$)  arc (135:225:\raycone);
\draw[thick,->,\colorboundcone] (b) -- ($(b) + (225:\lengthcone)$);
\draw[thick,->,\colorboundcone] (b) -- ($(b) + (135:\lengthcone)$);

\draw[thick,->,\colorboundcone]  ($ (c)!0.5!(b) $) -- ($(c)!0.5!(b) + (135:\lengthcone)$);

\fill[\colorcone] ($(c) + (135:\raycone)$) -- (c) -- ($(c)+ (45:\raycone)$)  arc (45:135:\raycone);
\draw[thick,->,\colorboundcone] (c) -- ($(c) + (45:\lengthcone)$);
\draw[thick,->,\colorboundcone] (c) -- ($(c) + (135:\lengthcone)$);

\draw[thick,->,\colorboundcone]  ($ (c)!0.5!(d) $) -- ($(c)!0.5!(d) + (45:\lengthcone)$);

\fill[\colorcone] ($(d) + (45:\lengthcone)$) -- (d) -- ($(d) + (-45:\lengthcone)$) arc (-45:45:\raycone);
\draw[thick,->,\colorboundcone] (d) -- ($(d) + (-45:\lengthcone)$);
\draw[thick,->,\colorboundcone] (d) -- ($(d) + (45:\lengthcone)$);

\draw[thick,->,\colorboundcone]  ($ (d)!0.5!(a) $) -- ($(d)!0.5!(a) + (-45:\lengthcone)$);

\draw[thin,dashed,->] (0,0) -- (1.25,0) node[anchor=west] {$y_1$};
\draw[thin,dashed,->] (0,0) -- (0,1.25) node[anchor=south] {$y_2$};
\draw[thin,dashed] (0,0) -- (-1.25,0);
\draw[thin,dashed] (0,0) -- (0,-1.25);

%% file: tikz/compfansupp_ball_infty.tex
\def\scalesub{\scalefibersfans}

\def\raycone{5.3}
\def\lengthcone{6.2}
\def\raybar{0.3}
\def\opacitycone{0.2}
\def\opacitypoly{0.3}
\def\rayreccone{2}
\def\lengthreccone{1.9}
\def\colorcone{green!80!white}
\def\colorboundcone{green!70!blue}
\def\colorpoly{orange!90!white}
\def\colorboundpoly{orange}
\def\colorcomplex{red}
\def\colorbar{pink}
\def\origin{0,0,0}

\def\R{3}

\def\ax{\R}      \def\ay{\R}
\def\bx{-\R}      \def\by{\R}
\def\cx{-\R}      \def\cy{-\R}
\def\dx{\R}      \def\dy{-\R}

\begin{subfigure}{.25\textwidth}
\centering

\begin{tikzpicture}[scale=\scalesub]

\coordinate (a) at (\ax,\ay);
\coordinate (b) at (\bx,\by);
\coordinate (c) at (\cx,\cy);
\coordinate (d) at (\dx,\dy);

\pgfmathsetmacro\anga{0};
\pgfmathsetmacro\angb{90};
\pgfmathsetmacro\angc{225};

\fill[fill=\colorpoly,opacity=\opacitypoly] (a) -- (b) -- (c) -- (d) --  (a);
\draw[\colorboundpoly,name path=polyhedron] (a)node{$\bullet$} -- (b)node{$\bullet$} -- (c)node{$\bullet$} -- (d)node{$\bullet$} -- (a);

\fill[\colorcone,opacity=\opacitycone] (\anga:\raycone) -- (0,0) --  (\angb:\raycone)  arc (\angb:\anga:\raycone);

\fill[\colorcone,opacity=\opacitycone] (\angb:\raycone) -- (0,0) -- (\angc:\raycone)  arc (\angc:\angb:\raycone);

\fill[\colorcone,opacity=\opacitycone] (\anga+360:\raycone) -- (0,0) -- (\angc:\raycone)  arc (\angc:\anga+360:\raycone);

\draw[thick,->,\colorboundcone,name path=halflinea] (0,0) -- (\anga:\lengthcone)node[black!60!green,anchor=west]{$\nhalfline{5}$};
\draw[thick,->,\colorboundcone,name path=halflineb] (0,0) -- (\angb:\lengthcone)node[black!60!green,anchor=south]{$\nhalfline{6}$};
\draw[thick,->,\colorboundcone,name path=halflinec] (0,0) -- (\angc:\lengthcone)node[black!60!green,anchor=north east]{$\nhalfline{3}$};

\draw[\colorcomplex,name intersections={of=polyhedron and halflinea}]
    (0,0)node{$\bullet$} -- (intersection-1)node{$\bullet$};

\draw[\colorcomplex,name intersections={of=polyhedron and halflineb}]
    (0,0)node{$\bullet$} -- (intersection-1)node{$\bullet$};

\draw[\colorcomplex,name intersections={of=polyhedron and halflinec}]
    (0,0)node{$\bullet$} -- (intersection-1)node{$\bullet$};

\draw[fill=\colorbar] (\R/2,\R/2) circle (\raybar);
\draw[fill=\colorbar] (-5*\R/9,2*\R/9) circle (\raybar);
\draw[fill=\colorbar] (2*\R/9,-5*\R/9) circle (\raybar);

\end{tikzpicture}
\caption{$\sigma=[-0.5,0]$}

\end{subfigure}
\begin{subfigure}{.25\textwidth}
\centering

\begin{tikzpicture}[scale=\scalesub]

\coordinate (a) at (\ax,\ay);
\coordinate (b) at (\bx,\by);
\coordinate (c) at (\cx,\cy);
\coordinate (d) at (\dx,\dy);

\pgfmathsetmacro\anga{0};
\pgfmathsetmacro\angb{90};
\pgfmathsetmacro\angc{135};
\pgfmathsetmacro\angd{225};
\pgfmathsetmacro\ange{315};

\fill[fill=\colorpoly,opacity=\opacitypoly] (a) -- (b) -- (c) -- (d) --  (a);
\draw[\colorboundpoly,name path=polyhedron] (a)node{$\bullet$} -- (b)node{$\bullet$} -- (c)node{$\bullet$} -- (d)node{$\bullet$} -- (a);

\fill[\colorcone,opacity=\opacitycone] (\anga:\raycone) -- (0,0) --  (\angb:\raycone)  arc (\angb:\anga:\raycone);

\fill[\colorcone,opacity=\opacitycone] (\angb:\raycone) -- (0,0) -- (\angc:\raycone)  arc (\angc:\angb:\raycone);

\fill[\colorcone,opacity=\opacitycone] (\angc:\raycone) -- (0,0) -- (\angd:\raycone)  arc (\angd:\angc:\raycone);

\fill[\colorcone,opacity=\opacitycone] (\angd:\raycone) -- (0,0) -- (\ange:\raycone)  arc (\ange:\angd:\raycone);

\fill[\colorcone,opacity=\opacitycone] (\anga+360:\raycone) -- (0,0) -- (\ange:\raycone)  arc (\ange:\anga+360:\raycone);

\draw[thick,->,\colorboundcone,name path=halflinea] (0,0) -- (\anga:\lengthcone)node[black!60!green,anchor=west]{$\nhalfline{5}$};
\draw[thick,->,\colorboundcone,name path=halflineb] (0,0) -- (\angb:\lengthcone)node[black!60!green,anchor=south]{$\nhalfline{6}$};
\draw[thick,->,\colorboundcone,name path=halflinec] (0,0) -- (\angc:\lengthcone)node[black!60!green,anchor=south east]{$\nhalfline{4}$};
\draw[thick,->,\colorboundcone,name path=halflined] (0,0) -- (\angd:\lengthcone)node[black!60!green,anchor=north east]{$\nhalfline{3}$};
\draw[thick,->,\colorboundcone,name path=halflinee] (0,0) -- (\ange:\lengthcone)node[black!60!green,anchor=north west]{$\nhalfline{2}$};

\draw[\colorcomplex,name intersections={of=polyhedron and halflinea}]
    (0,0)node{$\bullet$} -- (intersection-1)node{$\bullet$};

\draw[\colorcomplex,name intersections={of=polyhedron and halflineb}]
    (0,0)node{$\bullet$} -- (intersection-1)node{$\bullet$};

\draw[\colorcomplex,name intersections={of=polyhedron and halflinec}]
    (0,0)node{$\bullet$} -- (intersection-1)node{$\bullet$};

\draw[\colorcomplex,name intersections={of=polyhedron and halflined}]
    (0,0)node{$\bullet$} -- (intersection-1)node{$\bullet$};

\draw[\colorcomplex,name intersections={of=polyhedron and halflinee}]
    (0,0)node{$\bullet$} -- (intersection-1)node{$\bullet$};

\draw[fill=\colorbar] (-\R/3,2*\R/3) circle (\raybar);
\draw[fill=\colorbar] (\R/2,\R/2) circle (\raybar);
\draw[fill=\colorbar] (2*\R/3,-\R/3) circle (\raybar);
\draw[fill=\colorbar] (0,-2*\R/3) circle (\raybar);
\draw[fill=\colorbar] (-2*\R/3,0) circle (\raybar);

\end{tikzpicture}
\caption{$\sigma=[0,0.5]$}
\end{subfigure}
\begin{subfigure}{.24\textwidth}
\centering

\begin{tikzpicture}[scale=\scalesub]

\coordinate (a) at (\ax,\ay);
\coordinate (b) at (\bx,\by);
\coordinate (c) at (\cx,\cy);
\coordinate (d) at (\dx,\dy);

\pgfmathsetmacro\anga{0};
\pgfmathsetmacro\angb{45};
\pgfmathsetmacro\angc{90};
\pgfmathsetmacro\angd{135};
\pgfmathsetmacro\ange{225};
\pgfmathsetmacro\angf{315};

\fill[fill=\colorpoly,opacity=\opacitypoly] (a) -- (b) -- (c) -- (d) --  (a);
\draw[\colorboundpoly,name path=polyhedron] (a)node{$\bullet$} -- (b)node{$\bullet$} -- (c)node{$\bullet$} -- (d)node{$\bullet$} -- (a);

\fill[\colorcone,opacity=\opacitycone] (\anga:\raycone) -- (0,0) --  (\angb:\raycone)  arc (\angb:\anga:\raycone);

\fill[\colorcone,opacity=\opacitycone] (\angb:\raycone) -- (0,0) -- (\angc:\raycone)  arc (\angc:\angb:\raycone);

\fill[\colorcone,opacity=\opacitycone] (\angc:\raycone) -- (0,0) -- (\angd:\raycone)  arc (\angd:\angc:\raycone);

\fill[\colorcone,opacity=\opacitycone] (\angd:\raycone) -- (0,0) -- (\ange:\raycone)  arc (\ange:\angd:\raycone);

\fill[\colorcone,opacity=\opacitycone] (\ange:\raycone) -- (0,0) -- (\angf:\raycone)  arc (\angf:\ange:\raycone);

\fill[\colorcone,opacity=\opacitycone] (\anga+360:\raycone) -- (0,0) -- (\angf:\raycone)  arc (\angf:\anga+360:\raycone);

\draw[thick,->,\colorboundcone,name path=halflinea] (0,0) -- (\anga:\lengthcone)node[black!60!green,anchor=west]{$\nhalfline{5}$};
\draw[thick,->,\colorboundcone,name path=halflineb] (0,0) -- (\angb:\lengthcone)node[black!60!green,anchor=south west]{$\nhalfline{1}$};
\draw[thick,->,\colorboundcone,name path=halflinec] (0,0) -- (\angc:\lengthcone)node[black!60!green,anchor=south]{$\nhalfline{6}$};
\draw[thick,->,\colorboundcone,name path=halflined] (0,0) -- (\angd:\lengthcone)node[black!60!green,anchor=south east]{$\nhalfline{4}$};
\draw[thick,->,\colorboundcone,name path=halflinee] (0,0) -- (\ange:\lengthcone)node[black!60!green,anchor=north east]{$\nhalfline{3}$};
\draw[thick,->,\colorboundcone,name path=halflinef] (0,0) -- (\angf:\lengthcone)node[black!60!green,anchor=north west]{$\nhalfline{2}$};

\draw[\colorcomplex,name intersections={of=polyhedron and halflinea}]
    (0,0)node{$\bullet$} -- (intersection-1)node{$\bullet$};
\draw[\colorcomplex,name intersections={of=polyhedron and halflineb}]
    (0,0)node{$\bullet$} -- (intersection-1)node{$\bullet$};
\draw[\colorcomplex,name intersections={of=polyhedron and halflinec}]
    (0,0)node{$\bullet$} -- (intersection-1)node{$\bullet$};
\draw[\colorcomplex,name intersections={of=polyhedron and halflined}]
    (0,0)node{$\bullet$} -- (intersection-1)node{$\bullet$};
\draw[\colorcomplex,name intersections={of=polyhedron and halflinee}]
    (0,0)node{$\bullet$} -- (intersection-1)node{$\bullet$};
\draw[\colorcomplex,name intersections={of=polyhedron and halflinef}]
    (0,0)node{$\bullet$} -- (intersection-1)node{$\bullet$};

\draw[fill=\colorbar] (-\R/3,2*\R/3) circle (\raybar);
\draw[fill=\colorbar] (\R/3,2*\R/3) circle (\raybar);
\draw[fill=\colorbar] (2*\R/3,\R/3) circle (\raybar);
\draw[fill=\colorbar] (2*\R/3,-\R/3) circle (\raybar);
\draw[fill=\colorbar] (0,-2*\R/3) circle (\raybar);
\draw[fill=\colorbar] (-2*\R/3,0) circle (\raybar);

\end{tikzpicture}
\caption{$\sigma=[0.5,1]$}
\end{subfigure}
\begin{subfigure}{.24\textwidth}
\centering

\begin{tikzpicture}[scale=\scalesub]

\coordinate (a) at (\ax,\ay);
\coordinate (b) at (\bx,\by);
\coordinate (c) at (\cx,\cy);
\coordinate (d) at (\dx,\dy);

\pgfmathsetmacro\anga{45};
\pgfmathsetmacro\angb{135};
\pgfmathsetmacro\angc{225};
\pgfmathsetmacro\angd{315};

\fill[fill=\colorpoly,opacity=\opacitypoly] (a) -- (b) -- (c) -- (d) --  (a);
\draw[\colorboundpoly,name path=polyhedron] (a)node{$\bullet$} -- (b)node{$\bullet$} -- (c)node{$\bullet$} -- (d)node{$\bullet$} -- (a);

\fill[\colorcone,opacity=\opacitycone] (\anga:\raycone) -- (0,0) --  (\angb:\raycone)  arc (\angb:\anga:\raycone);

\fill[\colorcone,opacity=\opacitycone] (\angb:\raycone) -- (0,0) -- (\angc:\raycone)  arc (\angc:\angb:\raycone);

\fill[\colorcone,opacity=\opacitycone] (\angc:\raycone) -- (0,0) -- (\angd:\raycone)  arc (\angd:\angc:\raycone);

\fill[\colorcone,opacity=\opacitycone] (\anga+360:\raycone) -- (0,0) -- (\angd:\raycone)  arc (\angd:\anga+360:\raycone);

\draw[thick,->,\colorboundcone,name path=halflinea] (0,0) -- (\anga:\lengthcone)node[black!60!green,anchor=south west]{$\nhalfline{1}$};
\draw[thick,->,\colorboundcone,name path=halflineb] (0,0) -- (\angb:\lengthcone)node[black!60!green,anchor=south east]{$\nhalfline{4}$};
\draw[thick,->,\colorboundcone,name path=halflinec] (0,0) -- (\angc:\lengthcone)node[black!60!green,anchor=north east]{$\nhalfline{3}$};
\draw[thick,->,\colorboundcone,name path=halflined] (0,0) -- (\angd:\lengthcone)node[black!60!green,anchor=north west]{$\nhalfline{2}$};

\draw[thick,->,white] (90:\lengthcone)node[anchor=south]{$N_6$};

\draw[\colorcomplex,name intersections={of=polyhedron and halflinea}]
    (0,0)node{$\bullet$} -- (intersection-1)node{$\bullet$};
\draw[\colorcomplex,name intersections={of=polyhedron and halflineb}]
    (0,0)node{$\bullet$} -- (intersection-1)node{$\bullet$};
\draw[\colorcomplex,name intersections={of=polyhedron and halflinec}]
    (0,0)node{$\bullet$} -- (intersection-1)node{$\bullet$};
\draw[\colorcomplex,name intersections={of=polyhedron and halflined}]
    (0,0)node{$\bullet$} -- (intersection-1)node{$\bullet$};

\draw[fill=\colorbar] (0,2*\R/3) circle (\raybar);
\draw[fill=\colorbar] (0,-2*\R/3) circle (\raybar);
\draw[fill=\colorbar] (2*\R/3,0) circle (\raybar);
\draw[fill=\colorbar] (-2*\R/3,0) circle (\raybar);

\end{tikzpicture}
\caption{$\sigma=[1,+\infty)$}
\end{subfigure}

\centering
\caption{Exact quantization illustrated. The normal fan $\cN_\sigma$ in green with $N_i=W_i^\top \RR^+$, $\va c$ is uniform on the support $Q=-Q=B_\infty(0,R)$ in light orange, the sets $W_i^\top \RR^+\cap Q$ in red. The polyhedral complex $\compfansupp_\sigma$ shown in red or orange. The quantized costs $\check{c}_N$ are determined by centroids (small circles in pink).}

%% file: tikz/graph_Q_example.tex
\def\ax{-0.5}
\def\bx{0}
\def\cx{0.5}
\def\dx{1}
\def\ex{2}

\draw[thin,->] (-1,0) -- (2,0) node[anchor=west] {$x$};
\draw[thin,->] (0,-1.2) -- (0,0.2) node[anchor=south] {$\ectg(x)$};

\draw (\ax,0) node[anchor=south]{\ax};
\draw (\bx,0) node[anchor=south west]{\bx};
\draw (\cx,0) node[anchor=south]{\cx};
\draw (\dx,0) node[anchor=south]{\dx};

\def\param{1}
\pgfmathsetmacro\ay{(-7-14*\ax)/(8*\param)};
\pgfmathsetmacro\by{(-7-14*\bx)/(8*\param)};
\pgfmathsetmacro\cy{(-7-6*\cx)/(8*\param)};
\pgfmathsetmacro\dy{(-2-\dx)/(2*\param)};
\pgfmathsetmacro\ey{(-3)/(2*\param)};

\draw[blue] (\ax,\ay) -- (\bx,\by)-- (\cx,\cy)-- (\dx,\dy) -- (\ex,\ey)node[anchor=west] {$\frac{\expvalpar^2  e^{-\expvalpar \mynorm{c}_1}}{4} dc$};

\draw[thin,dashed] (\ax,0) -- (\ax,\ay);
\draw[thin,dashed] (\bx,0) -- (\bx,\by);
\draw[thin,dashed] (\cx,0) -- (\cx,\cy);
\draw[thin,dashed] (\dx,0) -- (\dx,\dy);

\def\param{1}
\pgfmathsetmacro\ay{(-7-14*\ax)*\param/24)};
\pgfmathsetmacro\by{(-7-14*\bx)*\param/24};
\pgfmathsetmacro\cy{(-7-6*\cx)*\param/24};
\pgfmathsetmacro\dy{(-2-\dx)*\param/6};
\pgfmathsetmacro\ey{-\param/2};

\draw[red] (\ax,\ay) -- (\bx,\by)-- (\cx,\cy)-- (\dx,\dy) -- (\ex,\ey)  node[anchor=south west] {uniform on norm 1 ball};

\draw[thin,dashed] (\ax,0) -- (\ax,\ay);
\draw[thin,dashed] (\bx,0) -- (\bx,\by);
\draw[thin,dashed] (\cx,0) -- (\cx,\cy);
\draw[thin,dashed] (\dx,0) -- (\dx,\dy);

\def\param{1}
\pgfmathsetmacro\ay{(-5-10*\ax)*\param/12)};
\pgfmathsetmacro\by{(-5-10*\bx)*\param)/12};
\pgfmathsetmacro\cy{(-5-4*\cx)*\param)/12};
\pgfmathsetmacro\dy{(-3-\dx)*\param/6};
\pgfmathsetmacro\ey{-2*\param/3};

\draw[green!60!black] (\ax,\ay) -- (\bx,\by)-- (\cx,\cy)-- (\dx,\dy) -- (\ex,\ey) node[anchor=north west] {uniform on norm $\infty$ ball};

\draw[thin,dashed] (\ax,0) -- (\ax,\ay);
\draw[thin,dashed] (\bx,0) -- (\bx,\by);
\draw[thin,dashed] (\cx,0) -- (\cx,\cy);
\draw[thin,dashed] (\dx,0) -- (\dx,\dy);

\def\param{1}
\pgfmathsetmacro\ay{-2*3.4142*(1+2*\ax)*\param/(3*2*3.1416)};
\pgfmathsetmacro\by{-2*3.4142*(1+2*\bx)*\param/(3*2*3.1416)};
\pgfmathsetmacro\cy{-2*(3.4142+2.8284*\cx)*\param/(3*2*3.1416)};
\pgfmathsetmacro\dy{-8*(1+0.4142*\dx)*\param/(3*2*3.1416)};
\pgfmathsetmacro\ey{-3.771*\param/(2*3.1416)};

\draw[orange!80!black] (\ax,\ay) -- (\bx,\by)-- (\cx,\cy)-- (\dx,\dy) -- (\ex,\ey) node[anchor=west] {uniform on norm 2 ball};

\draw[thin,dashed] (\ax,0) -- (\ax,\ay);
\draw[thin,dashed] (\bx,0) -- (\bx,\by);
\draw[thin,dashed] (\cx,0) -- (\cx,\cy);
\draw[thin,dashed] (\dx,0) -- (\dx,\dy);

\def\param{1}

\pgfmathsetmacro\ay{-3.4142*(1+2*\ax)*\param/(2*1.4142*1.7725)};
\pgfmathsetmacro\by{-3.4142*(1+2*\bx)*\param/(2*1.4142*1.7725)};
\pgfmathsetmacro\cy{-(3.4142+2.8284*\cx)*\param/(2*1.4142*1.7725)};
\pgfmathsetmacro\dy{-2*(1+0.4142*\dx)*\param/(1.4142*1.7725)};
\pgfmathsetmacro\ey{-2*\param/1.7725};

\draw[violet!90!black] (\ax,\ay) -- (\bx,\by)-- (\cx,\cy)-- (\dx,\dy) -- (\ex,\ey) node[anchor=west] {$\frac{e^{-\frac{\mynorm{c}_2^2}{2\gamma^2}}}{2\pi \gamma^2}dc$};

\draw[thin,dashed] (\ax,0) -- (\ax,\ay);
\draw[thin,dashed] (\bx,0) -- (\bx,\by);
\draw[thin,dashed] (\cx,0) -- (\cx,\cy);
\draw[thin,dashed] (\dx,0) -- (\dx,\dy);

%% file: table/affine_values_ectg.tex







\begin{table}[h!]
\centering
\begin{tabular}{c|c|c|c|c}
\footnotesize
$d\PP(c)$
& $-\frac{1}{2}\leq x \leq 0$
& $0\leq x \leq\frac{1}{2}$
&  $\frac{1}{2} \leq x \leq1$
& $1 \leq x $ \\ \hline

$\frac{\indi{\mynorm{c}_1 \leq R}}{2R^2} dc$
& $\frac{-7R}{24} (1+2x)$
& $\frac{-R}{24} (7+6x)$
& $\frac{-R}{6} (2+x)$
& $\frac{-R}{2}$\\ \hline

$\frac{\expvalpar^2  e^{-\expvalpar \mynorm{c}_1}}{4} dc$
& $\frac{-7}{8 \expvalpar} (1+2x)$
& $\frac{-1}{8 \expvalpar}(7+6x)$
& $\frac{-1}{2 \expvalpar}(2+x)$
& $\frac{-3}{2 \expvalpar}$\\ \hline

$\frac{\indi{\mynorm{c}_\infty \leq R}}{4R^2} dc$ 
& $\frac{-R}{12} (5+10x)$
& $\frac{-R}{12} (5+4x)$ 
& $\frac{-R}{6} (3+x)$
& $\frac{-2R}{3}$\\ \hline

$\frac{e^{-\mynorm{c}_2^2/2\gamma^2}}{2\pi \gamma^2}dc$
& $\frac{-\gamma(2+\sqrt2)(1+2x)}{2\sqrt{2\pi}}$
& $\frac{-\gamma(2+\sqrt2 +2\sqrt2 x)}{2\sqrt{2\pi}}$
& $\frac{-2\gamma(1+(-1+\sqrt2)x)}{\sqrt{2\pi}}$
& $-\frac{2}{\sqrt{\pi}}\gamma$\\ \hline 

$\frac{\indi{\mynorm{c}_2 \leq R}}{\pi R^2} dc$
& $\frac{-R(2+\sqrt2)(1+2x)}{3\pi}$
& $\frac{-R(2+\sqrt2 +2\sqrt2 x)}{3\pi}$
& $\frac{-4R(1+(-1+\sqrt2)x)}{3\pi}$
& $-\frac{4\sqrt2 R}{3\pi}$
\\ 
\end{tabular}
\caption{Different values of $\ectg(x)$ for different distribution of the cost $\va c$}
\label{tab:values_ectg}
\end{table}

%% file: sections/complexity.tex
Hanasusanto, Kuhn and Wiesemann
showed in~\cite{hanasusanto2016comment} that 
2-stage stochastic programming is $\sharp$P-hard,
by reducing the computation of the volume of a polytope to the 
resolution of a 2-stage stochastic program.
Nevertheless, we show that for a fixed dimension of the recourse space,
2-stage programming is polynomial.
Therefore, the status of 2-stage programming seems
somehow comparable to the one of the computation
of the volume of a polytope -- which is also both $\sharp$P-hard
and polynomial
when the dimension is fixed (see for example \cite[3.1.1]{gritzmann1994complexityii}).
We also give a similar result for multistage stochastic linear programming.

We now give a summary of our method.
A naive approach would be to use directly
the exact quantization result \cref{thm:quantization_2stage},
for every $x$.
However, even in the two stage case, the latter yields a linear program
  of an exponential size when only the recourse dimension $m$ is fixed.
  Indeed, the size of the quantized linear program, $(2SLP)$
  is polynomial only when {\em both} $n$ and $m$ are fixed.
  Indeed, $\bigwedge_{\sigma \in \chcmplx{P}{\prjx}} -\cN_\sigma$ can have, by McMullen's and Stanley's upper bound theorems (\cite{mcmullen1970maximum,stanley1975upper}),
  an exponential size in $n$ and $m$, and these bounds are tight. 
  Hence, to handle the case in which only the recourse dimension $m$ is fixed, we need additional ideas.
   We use the quantization result \cref{thm:quantization_2stage} only for a {\em fixed} $x$, observing that when $m$ is fixed, $\cN(P_x)$ has a polynomial size.  We thus have a polynomial time
   oracle that gives the values $\ectg(x)$ by \cref{thm:quantization_2stage} and a subgradient $g \in \partial \ectg(x)$ by \cref{cor:subgradient_2stage}.
Then, we rely on the theory of linear programming with oracle~\cite{grotschel2012geometric}, working in the \emph{Turing model} of computation (a.k.a. \emph{bit model}). In particular, all the computations are carried out with rational numbers. We now provide the needed details
of the proof.

\subsection{Multistage programming with exact oracles}

\label{ssec:exact_rational}
Recall that a polyhedron can be given in two manners.
The ``$H$-representation'' provides an external description
of the polyhedron, as the intersection of finitely
many half-spaces. The ``$V$-representation'' provides
an internal representation, writing the polyhedron
as a Minkowski sum of a polytope (given as the convex
hull of finitely many points) and of a polyhedral
cone (generated by finitely many vectors).  

We say that a polyhedron is {\em rational} if 
the inequalities in its $H$-representation are rational or, equivalently, the generators of its $V$-representation have rational coefficients.
We shall say that a (convex) polyhedral function $V$ is {\em rational} if its epigraph
is a rational polyhedron.

Recall that, in the Turing model, the {\em size} (or encoding length see \cite[1.3]{grotschel2012geometric}) of an integer $k \in \ZZ$ is $\enclgth{k} \defegal 1 + \lceil\log_2(|k|+1)\rceil$;
the size of a rational $r=\frac{p}{q} \in \QQ$ with $p$ and $q$ coprime integers, is $\enclgth{r} \defegal \enclgth{p} +\enclgth{q}$.
The size of a rational matrix or a vector, still denoted by $\enclgth{\cdot}$,
is the sum of the sizes of its entries.
The size of an inequality $\alpha^\top x\leq \beta$ is $\enclgth \alpha + \enclgth \beta$.
The size of a $H$-representation of a polyhedron is the sum of the sizes of its inequalities and the size of a $V$-representation of a polyhedron is the sum of the sizes of its generators.

If the dimension of the ambient space is {\em fixed},
one can pass from one representation to the other one
in {\em polynomial time}.
Indeed, the double description algorithm allows one to get
a $V$-representation from a $H$-representation, 
see the discussion at the end of section 3.1 in~\cite{fukuda1995double}, and use McMullen's upper bound theorem (\cite{mcmullen1970maximum} and \cite[6.2.4]{grotschel2012geometric}) to show that the computation
time 
is polynomially bounded in the size of the $H$-representation.
A fortiori, the size of the $V$-representation is polynomially
bounded in the size of the $H$-representation.
Dually, the same method allows one to obtain a $H$-representation from a $V$-representation.
Hence, in the sequel, we shall use the term {\em size} of a polyhedron
for the size of a $V$ or $H$-representation: when dealing
with polynomial-time complexity results in fixed dimension,
whichever representation is used is irrelevant. In particular,
we define the {\em size} $\enclgth{N}$ of a rational cone $N$
as the size of a $H$ or $V$ representation of $N$.

We first observe that the size of the scenario tree arising in the
exact quantization result becomes polynomial when suitable dimensions
are fixed.

\begin{proposition}\label{prop:size-ptime}
Let $t\in \{2,\dots,\horizon\}$, and suppose that the dimensions $n_t,\dots,n_{\horizon}$ and the cardinals $\card (\supp \va \xi_t)$, $\cdots$,$\card (\supp \va \xi_{\horizon})$ are fixed. Let $\cT$ be the scenario tree constructed in~\Cref{cor:quantization_scenario_tree}. Then, the subtree of $\cT$ rooted
at an arbitrary node of depth $t$ can be computed in polynomial in $\sum_{s=t}^{\horizon} \sum_{\xi\in\supp(\va\xi_s)} \enclgth{\xi}$.
\end{proposition}
\begin{proof}
Recall that
the number of chambers of a chamber complex is polynomial when both dimensions are fixed by \cite[3.9]{verdoolaege2005computation}.
Thus, we can compute recursively the (maximal) chambers of the complexes $\cP_t$ defined in \cref{thm:polyhedral_multistage} thanks to the algorithm in \cite[3.2]{clauss1998parametric} in polynomial time.
We then can compute in polynomial time the fans $\cN_t$ defined in \cref{thm:multistage_quantization}.
\end{proof}

We recall the theory of linear programming with oracle applies to the class of ``well described'' polyhedra  which are rational polyhedra with an apriori bound on the bit-sizes of the inequalities defining their facets, we refer the reader to~\cite{grotschel2012geometric}
for a more detailed discussion of the notions (oracles) and results
used here.

  \begin{definition}[first-order oracle]
  Let $f$ be a rational polyhedral function. We say that $f$ admits a polynomial time (exact) \emph{first-order oracle}, if there exists an oracle that takes as input a vector $x$ and either returns a hyperplane separating $x$ from $\dom(f)$ if $x \notin \dom(f)$ or returns $f(x)$ and $g \in \partial\ectg(x)$ if $x \in \dom(f)$, in polynomial time in $\enclgth{x}$.
    \label{defi:first_order_oracle}    
  \end{definition}

\begin{lemma}
\label{lem:first_order_to_polynomial_complexity}
   Let $Q \subset \RR^d$ be a polyhedron, $c \in \RR^d$ a cost vector and $f$ be a polyhedral function given by a first-order oracle.
   Futhermore, assume $\epi(f)$ and $Q$ are well described.
   Then, the problem $\min_{x\in Q} \; c^\top x +f(x)$ can be solved in oracle-polynomial time in $\enclgth{c}+\enclgth{\epi(f)}+\enclgth{Q}$.
\end{lemma}

\begin{proof}
  The case where $\dom(f)=\RR^d$ is tackled in
  Theorem 6.5.19 in \cite{grotschel2012geometric}.
  If $f$ has a general domain, we can write $f=\tilde f+\tropindi{\dom f}$
  where $\tilde f$ is a polyhedral function with a well described
  epigraph and such that $\dom \tilde f=\RR^d$.
  Then, noting that $\epi(f)=\epi(\tilde f) \cap \dom(f)\times \RR$, we can adapt the proof of the latter theorem, using Exercise 6.5.18(a) of \cite{grotschel2012geometric}.
\end{proof}

We do not require the distribution of the cost $\va c$ to be
described extensively. We only need to assume the
existence of the following oracle.

  \begin{definition}[cone-valuation oracle]
  Let $\va c\in L(\Omega,\cA,\PP,\RR^m)$ be an integrable
   cost distribution
   such that,
   for every \emph{rational} cone $N$,  the quantized probability $\check p_N$ and  quantized cost $\check c_N$ are 
    \emph{rational}.
    We say that $\va c$ admits a polynomial time (exact) \emph{cone-valuation oracle}, if there exists an oracle which takes as input
    a rational polyhedral cone $N$ and returns $\check p_N$ and $\check c_N$ in \emph{polynomial time} in $\enclgth{N}$. 
    \label{defi:oracle_cost}    
  \end{definition}

\begin{theorem}[Cone valuation to first-order oracle]
\label{thm:cone_valuation_to_first_order}
Consider the value functions of problem \eqref{eq:def_mslp} defined in \cref{eq:def_ectg_multi} . 
Assume that $\horizon,n_2,\dots,n_{\horizon} $,  $\card (\supp \va \xi_2)$, $\cdots$,$\card (\supp \va \xi_{\horizon})$ are fixed integers, and that 
$(\va c_t, \va \xi_t)_{2 \leq t \leq \horizon}$ 
satisfies \cref{as:supp_vac_multi}.
Assume in addition that, every vector $\xi \in \supp(\va \xi_t)$ has rational
entries and that the probabilities $p_{t,\xi}:=\bprob{\va \xi_t=\xi}$
are rational numbers.
Assume finally that every random variable $\va c_t$ conditionally to $\{\va \xi_t = \xi \}$, denoted by $\va c_{t,\xi}$, admits a polynomial-time cone-valuation oracle (see \cref{defi:oracle_cost}).

Then, for all $t\geq 2$, $\ectg_t$ admits a polynomial time first-order oracle.

\end{theorem}

\begin{proof}
We start with the 2-stage case with deterministic constraints. We recall our notation $\ectg(x):=\besp{\min_{y \in \RR^m}\va c^\top y +\tropindi{Ay+ Bx\leq b}}$.
Let $x \in \RR^n$ be an input vector. 
We first check if $x\in \prj(P)=\dom(V)$.
By solving the dual of $\min_{y\in \RR^q} \{ 0 \, | \, Ay \leq b-Bx\}$,
we either 
find an unbounded ray generated by $\lambda \in \RR^q$ such that $\lambda\geq 0$, $\lambda^\top A=0$ and $\lambda^\top (b-Bx)<0$
or a $y\in \RR^m$ such that $Ay\leq b-Bx$, so that $x\in \prj(P)$.
In the former case we have $x\notin \prj(P)$, and we get a cut
$\{ x'\in \RR^{n} \, |\, \lambda^\top Bx' =\frac{\lambda^\top b+\lambda^\top Ax}{2}\}$, separating $\prj(P)=\dom(\ectg)$ from $x$.

So, we  now assume that $x\in \prj(P)$, i.e., $V(x)<+\infty$.
We next show that we can compute $\ectg(x)$ and a subgradient $\alpha \in \partial \ectg(x)$ in polynomial time. 
Indeed, the McMullen upper-bound theorem \cite{mcmullen1970maximum}, in its dual version, guarantees that a polytope of dimension $m$ with $f$ facets has $O(f^{\lfloor m/2 \rfloor})$ faces, see~\cite{Sei95}.
  Since the number of cones in $\cN(P_x)$ is equal to the number of faces of $P_x$ which is polynomially bounded in the number of constraints $q \leq \enclgth{\xi}$,   $\card \cN(P_x)$ is polynomial in $\enclgth{\xi}$.
    Thus, since $\va c$ is given by a cone valuation oracle, we can compute in polynomial time the collection of all quantized costs and probabilities
    $\check c_N$ and $\check p_N$, indexed by $N \in -\cN(P_x)$.
    Then, by \cref{thm:quantization_2stage}, we can compute $\ectg(x)$ by solving a linear program for each cone $N\in - \cN(P_x)$.
    Similarly, \cref{cor:subgradient_2stage}, allows us to
    compute a subgradient $\alpha \in
    \partial\ectg(x)$ using the same linear programs.
  All these operations take a polynomial time.

  The case of finitely supported
stochastic constraints reduces to the case of deterministic
constraints dealt with above,
using $\dom(\ectg)=\cap_{\xi \in \supp\va \xi} \prj(P(\xi))$ and
$\ectg(x)=\sum_{\xi \in \supp \va \xi} p_\xi \tilde \ectg(x|\xi)$
where $\tilde \ectg(x|\xi):=\besp{\ctg(x,\va c,\va \xi)\; |\; \va \xi =\xi}$.

We finally deal with the multistage case in a similar way, using the quantization result \Cref{cor:quantization_scenario_tree} in extensive form.
Applying~\Cref{prop:size-ptime}, the quantized costs
and probabilities arising there can be computed by a polynomial
number of calls to the cone-valuation oracle. This provides
a first order oracle for the expected cost-to-go function $\ectg_t$.

\end{proof}

We now refine the definition of cone-valuation
oracle, to take into account situations in which the distribution
of the random cost $\va{c}$ is specified by a parametric model.
We shall say that such a distribution admits
a polynomial-time {\em parametric cone-valuation oracle}
if there is an oracle that takes as input the parameters
of the distribution, together with a rational cone $N$,
and outputs the quantized probability $\check{p}_N$ and cost
$\check{c}_N$. Especially, we consider the following
situations:
\begin{itemize}
      \item[1.] {\em Deterministic distribution} equal to a rational cost $c$.
      We set $\enclgth{\va c}:=\enclgth{c}$
      \item[2.] {\em Exponential distribution on a rational cone} $K$ with rational parameter $\theta$. We set $\enclgth{\va c}:=\enclgth{K}+ \enclgth{\theta}$
      \item[3.] {\em Uniform distribution on a rational polyhedron $Q$}
      such that 
$\aff(Q)=\{y \in \RR^m \;|\; \forall j \in J\subset [m], y_j=q_j \in \QQ \}$ where $J$ is a subset of $[m]$ and $q_j$ are rational numbers (in particular, $Q$ is full dimensional when $J=\emptyset$). We set: $\enclgth{\va c}=\enclgth{Q}$
      \item[4.] {\em Mixtures of the above distributions}, i.e., convex combination with rational coefficients $(\lambda^k)_{k \in [l]}$ of distributions of random variables $(\va c_k)_{k \in [l]}$ satisfying 1. 2. or 3. Then, we set $\enclgth{\va c}=\sum_{k=1}^l\enclgth{\va c_k}+\enclgth{\lambda_k}$.
\end{itemize}
  \begin{theorem}
  Assume that the dimension $m$ is fixed, and that $\va c$ is
  distributed according to any of the above laws (deterministic, exponential,
  uniform, or mixture). Then, 
  the random cost $\va{c}$ admits a polynomial-time parametric cone-valuation oracle.
    \label{thm:usual_distr_cone_valuation_oracle}
  \end{theorem}
\begin{proof}
1. {\em Case of a deterministic distribution}.
We first check whether $c \in \relint(N)$,
which can be done in polynomial time, see section~6.5 of~\cite{grotschel2012geometric}.
Then,  if $c \in \relint(N)$, we set $\check c_N=c$ and $\check p_N=1$ otherwise $\check c_N=0$ and $\check p_N=0$.

2. {\em Case of an exponential distribution}.
Since the dimension is fixed, 
  for every polyhedron $R$, we can triangulate $R \cap \supp(\va c)$ and partition it into (relatively open) simplices and simplicial cones $(S_k)_{k\in[l]}$,
and by Stanley upper bound theorem, the size $l$ of the triangulation is polynomial in $\enclgth{R}$.
  By using the Brion formula in \cref{tab:synthesis_formulas_cost},
  we compute in polynomial time $\check p_R = \sum_{k=1}^l \check p_{S_k}$ and $\check c_R =  \sum_{k=1}^l \check p_{S_k} \check c_{S_k} / \check p_R$ if $\check p_R= 0$ and $\check c_R=0$ otherwise.

3. {\em Case of a uniform distribution}.
After triangulating (as in the case of an exponential distribution),
we may suppose that the support of the distribution
is a simplex $S$, so that $Q=S$.
If this simplex $S$ is full dimensional, then its volume
is given by a determinantal expression,
and so, it is rational (see e.g.~\cite{gritzmann1994complexityii} 3.1).
Then, the formulas of \Cref{tab:synthesis_formulas_cost} yield the result.
If this simplex is not full dimensional,
we have $\aff(S)=\{y \in \RR^m \;|\; \forall j \in J, y_j=q_j\}$,
a similar formula holds, ignoring the coordinates of $y$ whose indices are in the set
$J$.

4. {\em Case of mixtures of distributions}. Trivial reduction
to the previous cases.
\end{proof}

\begin{remark}
The conclusion of~\Cref{thm:usual_distr_cone_valuation_oracle}
does not carry over to the uniform distribution on
a general polytope of dimension $k<n$. The condition
that $\aff(Q)=\{y \in \RR^m \;|\; \forall j \in J, y_j=q_j\}$
ensures that the orthogonal projection on $\aff(Q)$ preserves
rationality, which entails that the $k$-dimensional
volume of $Q$ is a rational number.
In general, this volume is obtained by applying the Cayley Menger determinant formula (see for example \cite[3.6.1]{gritzmann1994complexityii}),
and it belongs to a quadratic extension of the field of rational numbers.
For example, if  $\Delta_d$ is the canonical simplex $\{\lambda\in \RR_+^{d+1}| \sum_{i =1}^{d+1}\lambda_i=1\}$ then $\vol(\Delta_d)=\frac{\sqrt{d+1}}{d!}$.

For the Gaussian distribution,
$\check c_S$ and $\check p_S$ can be determined in terms
of solid angles (see \cite{ribando2006measuring})  arising in~\Cref{tab:synthesis_formulas_cost}.
These coefficients are generally involving the number $\pi$ and Euler's $\Gamma$ function, and thus they are irrational.
\end{remark}

\begin{cor}[MSLP is polynomial for fixed dimensions]
Consider the problem \cref{eq:def_mslp} . 
Assume that $\horizon,n_2,\dots,n_{\horizon} $,  $\card (\supp \va \xi_2)$, $\cdots$,$\card (\supp \va \xi_{\horizon})$ are fixed integers, that 
$(\va c_t, \va \xi_t)_{2 \leq t \leq \horizon}$ 
satisfies \cref{as:supp_vac_multi}.
Suppose in addition that, for all $\xi \in \supp(\va \xi_t)$, $p_{t,\xi}:=\bprob{\va \xi_t=\xi}$ and $\xi$ are rational and that the random variable $\va c_t$ conditionally to $\{\va \xi_t = \xi \}$, denoted by $\va c_{t,\xi}$, is of the type considered in \cref{thm:usual_distr_cone_valuation_oracle}.

Then, Problem \eqref{eq:def_mslp} can be solved in a time that is polynomial in the input size $\enclgth{c_1}+\enclgth{\xi_1} +\sum_{t=2}^{\horizon} \sum_{\xi \in \supp(\va \xi_t)} (\enclgth{\va c_{t,\xi}} + \enclgth{\xi}+ \enclgth{p_{t,\xi}})$.
\end{cor}

\begin{proof}
We first show by backward induction that the epigraph $\epi(\ectg_2)$ is well described. 
The dynamic programming equation \cref{eq:def_ectg_multi}
allows us to compute a $H$-representation of $\epi(\ectg_t)$
from a $H$-representation of $\epi(\ectg_{t+1})$.
Indeed, by \cref{thm:multistage_quantization}, we have 
\begin{align*}
\ectg_t(x_{t-1}) &= 
    \sum_{\xi \in \supp(\va \xi_t)} p_{t,\xi}
   \sum_{N \in \cN_{t,\xi}} \check p_{t,N|\xi} \; \min_{x_t \in \RR^{n_t}} Q_{t,N|\xi}(x_t,x_{t-1})\enspace,\textrm{ with}\\
 Q_{t,N|\xi}(x_t,x_{t-1})&:=\check{c}_{t,N|\xi}^\top x_t  + \ectg_{t+1}(x_t)+
    \tropindi{ (x_t,x_{t-1})\in P_t(\xi)}
    \enspace .
    \end{align*}
We then have
\begin{align*} \epi(Q_{t,N|\xi})&=\bp{\epi(x_t\mapsto\check{c}_{t,N|\xi}^\top x_t)+\epi(\ectg_{t+1})} \; \cap \; (P_t(\xi)\times \RR)
\\
 \epi (\ectg_t) &= 
\sum_{\xi \in \supp(\va \xi_t)} p_{t,\xi}
   \sum_{N \in \cN_{t,\xi}} \check p_{t,N|\xi} \;\prj^{x_{t-1},x_t,z}_{x_{t-1},z}\!\big(\epi(Q_{t,N|\xi})\big)
   \enspace ,
   \end{align*}
   recalling that $\prj^{x_{t-1},x_t,z}_{x_{t-1},z}$ denotes the projection
   mapping $(x_{t-1},x_{t},z)\mapsto (x_{t-1},z)$. Well described
   polyhedra are stable under the operations of projection,
   intersection, and Minkowski sum, see in particular~\cite[6.5.18]{grotschel2012geometric}. It follows that
$\epi (\ectg_t)$ is well described.
Then, the corollary follows from~\Cref{lem:first_order_to_polynomial_complexity}, \Cref{thm:cone_valuation_to_first_order} and ~\Cref{thm:usual_distr_cone_valuation_oracle}.
\end{proof}

\subsection{Multistage programming with inexact oracles}
We finally consider the situation in which the law
of the cost distribution is only known approximately.
Hence, we relax the notion of cone-valuation oracle,
as follows.

  \begin{definition}[Weak cone-valuation oracle]
  Let $\va c\in L(\Omega,\cA,\PP,\RR^m)$ be an integrable
   cost distribution.
    We say that $\va c$ admits a polynomial time \emph{weak cone-valuation oracle}, if there exists an oracle which takes as input
    a rational polyhedral cone $N$ together with a rational number
    $\varepsilon>0$, and returns a rational number $\tilde{p}_N$ and
    a rational vector $\tilde{c}_N$ such that $|\tilde{p}_N-\check p_N|\leq \varepsilon$ and $\|\tilde{c}_N-\check c_N\|\leq \varepsilon$, in a time
    that is \emph{polynomial} in $\enclgth{N} + \enclgth{\varepsilon}$.
    \label{defi:oracle_cost_weak}
  \end{definition}

\begin{definition}[Weak first-order oracle]
  Let $f$ be a rational polyhedral function. We say that $f$ admits a polynomial time \emph{weak first-order oracle}, if there exists an oracle that takes as input a vector $x$ and either returns a hyperplane separating $x$ from $\dom(f)$ if $x \notin \dom(f)$ or returns a scalar $\tilde f$ and a vector $\tilde g$ such that $|\tilde f -f(x)| \leq \varepsilon$ and $d\bp{\tilde g,\partial f(x)}\leq \varepsilon$ if $x \in \dom(f)$, in a time which is polynomial in $\enclgth{x}+\enclgth{\varepsilon}$.
    \label{defi:first_order_oracle_weak}   
\end{definition}
\begin{remark}
In our definition of weak first order oracle, we require that feasibility
($x\in \dom(f)$)
be tested exactly, whereas the value and a subgradient of the function
are only given approximately. This is suitable to the present
setting, in which the main
difficulty resides in the approximation of the function (which may take
irrational values for relevant cost distributions).
\end{remark}
We now rely on the theory of linear programming with
weak separation oracles developed in~\cite{grotschel2012geometric}.
Let $C\subset \RR^d$ be convex set, for $\varepsilon>0$,
let $S(C,\varepsilon):=\{x\in \RR^d\mid \|x-y\|\leq \varepsilon\}$
and $S(C,-\varepsilon):=\{x\in \RR^d\mid B(x,\varepsilon)\subset C\}$
where $B(x,\varepsilon)$ denotes the Euclidean ball centered
at $x$ of radius $\varepsilon$. 
A {\em weak separation oracle}
for a convex set $C \subset \RR^d$ takes as argument a vector $x \in \RR^d$
and a rational number $\varepsilon>0$,
and either asserts that $x \in S(C,\varepsilon)$ or returns
a rational vector $\gamma\in \RR^d$, of norm one, and a rational scalar $\delta$, such that $\gamma^\top y \leq\gamma^\top x+\varepsilon$
for all $y\in S(C,-\varepsilon)$. 

\begin{theorem}[Weak cone valuation to weak first-order oracle]
\label{thm:cone_valuation_to_first_order_weak}
Consider the value functions of problem \eqref{eq:def_mslp} defined in \cref{eq:def_ectg_multi} . 
Assume that $\horizon,n_2,\dots,n_{\horizon} $,  $\card (\supp \va \xi_2)$, $\cdots$,$\card (\supp \va \xi_{\horizon})$ are fixed integers, and that 
$(\va c_t, \va \xi_t)_{2 \leq t \leq \horizon}$ 
satisfies \cref{as:supp_vac_multi}.
Assume in addition that, every vector $\xi \in \supp(\va \xi_t)$ has rational
entries and that the probabilities $p_{t,\xi}:=\bprob{\va \xi_t=\xi}$
are rational numbers. 
Assume finally that the diameters
of $\dom \ectg_t$, for $t\geq 2$, are bounded
by a rational constant $R$, and that 
every random variable $\va c_t$ conditionally to $\{\va \xi_t = \xi \}$, denoted by $\va c_{t,\xi}$, admits a polynomial-time weak cone-valuation oracle (see \cref{defi:oracle_cost}).

Then, for all $t\geq 2$, $\ectg_t$ admits a polynomial time weak first-order oracle.

\end{theorem}
\begin{proof}
The proof is similar to the one of~\Cref{thm:cone_valuation_to_first_order}.
The main difference is that we need an apriori bound $R$
on the diameter of $\dom \ectg_t$, so that if $d(\tilde{g},\partial \ectg_t(x))\leq \varepsilon$, then, using Cauchy-Schwarz inequality, $\ectg_t(y)-\ectg_t(x)\geq \tilde{g}\cdot (y-x) - \varepsilon R$ holds for all $y\in \dom\ectg_t$. Together with
and approximation of $\ectg_t(x)$, this allows us to get a weak separation oracle for the epigraph of $\ectg_t$.
\end{proof}

    \begin{cor}[Approximate (MSLP) is polynomial-time for fixed recourse dimension $m$]
Consider Problem \eqref{eq:def_mslp}.
Assume that $\horizon,n_2,\dots,n_{\horizon} $,  $\card (\supp \va \xi_2)$, $\cdots$,$\card (\supp \va \xi_{\horizon})$ are fixed integers 
Assume finally that the diameters
of $\dom \ectg_t$, for $t\geq 2$, are bounded
by a rational constant $R$, and that for all $\xi \in \supp(\va \xi_t)$, the random variable $\va c_t$ conditionally to $\{\va \xi_t = \xi \}$, denoted by $\va c_{t,\xi}$, admits a polynomial-time weak cone-valuation oracle.

Then, there exists an algorithm that either asserts that Problem \cref{eq:def_mslp} is infeasible or find a feasible solution $x^*$ whose cost
does not exceed the cost of an optimal solution by more than $\varepsilon$, in polynomial-time in $\enclgth{c_1}+\enclgth{\xi_1} +\sum_{t=2}^{\horizon} \sum_{\xi \in \supp(\va \xi_t)} (\enclgth{\va c_{t,\xi}} + \enclgth{\xi}+ \enclgth{p_{t,\xi}})+\enclgth{R}$.
\label{thm:complexity_2stage_m_fixed_oracle_weak}
\end{cor}
\begin{proof}
This follows from~\Cref{thm:cone_valuation_to_first_order_weak}, using
the result analogous to~\Cref{lem:first_order_to_polynomial_complexity}
for weak separation oracles, see~\cite[6.5.19]{grotschel2012geometric}.
\end{proof}
Finally, we show that every absolutely continuous cost distribution,
with a suitable density function,
admits a polynomial-time weak cone-valuation oracle.
\begin{definition}
We shall say that a density function $f:\RR^{n}\to \RR_{+}$
is {\em combinatorially tight} if:
\begin{enumerate}
\item there is a polynomial
time algorithm which, given a rational number $\varepsilon>0$,
returns a rational number $r>0$ such that $\int_{\|x\|>r} f(x)dx \leq \varepsilon$.
\item there is a polynomial time algorithm, which given
a rational vector $x\in \RR^n$, and a rational number $\varepsilon>0$,
returns an $\varepsilon$ approximation of $f(x)$.
\end{enumerate}\label{def:tight}
\end{definition}
The terminology is inspired by the notion of tightness from measure theory
(analogous to condition~1 in~\Cref{def:tight}).

We shall need a classical result on the numerical approximation
of multidimensional integrals, which can be found in~\cite{davis1984methods}.
The {\em total variation in the sense of Hardy and Krause}, $\|f\|_{\textrm{BVHK}}$, of a 
function $f$ on a $n$ dimensional hypercube
is defined in \cite[Def. ~p.352]{davis1984methods}). In particular,
if $f$ is of regularity class $\mathcal{C}^n$, $\|f\|_{\textrm{BVHK}}$
is finite. The error made when approximating
the integral of a function of $n$ variables by its Riemann sum taken on a regular grid with $k$ points
is bounded by $(n\|f\|_{\textrm{BVHK}})/k^{1/n}$, see the theorem on p~352 of~\cite{davis1984methods}.

\begin{proposition}\label{prop:yeswecan}
Suppose that a cost distribution $\va c$ admits a density function $f:\RR^{n}\to \RR_{+}$, that
is such that the function
$(1+\|\cdot\|)f$ is combinatorially tight
and that it has a finite total variation in the sense of Hardy and Krause,
bounded by an a priori constant.
Suppose
that the dimension $n$ is fixed.
Then,
$\va c$ admits a polynomial-time weak cone valuation oracle.
\end{proposition}

\begin{proof}
Given a rational cone $N$, we need to approximate
the integrals $\int_N f(c) dc$ and $\int_N cf(c) dc$, up
to the precision $\varepsilon$. Using the tightness condition,
it suffices to approximate the integrals of the same functions
restricted to the domain $N_r:= N\cap B_\infty(0,r)$, where $B_\infty(0,r)$
denotes the sup-norm ball of radius $r$, and the encoding length
of $r$ is polynomially bounded in the encoding length of $\varepsilon$.
We only discuss the approximation of $\int_{N_r} cf(c)dc$
(the case of $\int_{N_r} f(c)dc$ being simpler).
We denote by $\tilde{c}_{N_r}$
the approximation of $\int_{N_r} cf(c)dc$ provided
by taking the Riemann sum of the function $c\mapsto cf(c)$ over
the grid $([-r,r))^n\cap ((r/M)\mathbb{Z})^n$, which has $(2M)^r$ points.
Then, setting $g:=(1+\|\cdot\|)f$,
it follows from the result~\cite[Th.~p~352]{davis1984methods} recalled above that
$\|\int_{N_r} cf(c)dc-\tilde{c}_{N_r}\| \leq n \|g\|_{\textrm{BVHK}}/(2M)$.
Hence, for a fixed dimension $n$, we can get an $\varepsilon$
approximation of $\int_N cf(c) dc$ in a time polynomial in the
encoding length of $\varepsilon$.
\end{proof}
\begin{remark}
\Cref{prop:yeswecan} and~\Cref{thm:complexity_2stage_m_fixed_oracle_weak}
entail that, under the previous fixed-parameter restrictions
(including dimensions of the recourse spaces),
the MSLP problem is polynomial-time approximately
solvable for a large class of cost distributions. This applies
in particular to distributions like Gaussians, which are combinatorially tight.
In this case, condition~1 of~\Cref{def:tight}, whereas
condition~2 follows from the result of~\cite{borwein},
implying that the exponential function, restricted to the interval
$(-\infty,0]$, can be approximated in polynomial time.
\end{remark}